\documentclass[a4paper,12pt]{article}
\usepackage{amsthm,amsmath,amssymb}
\usepackage{epsfig}
\usepackage{geometry}
\usepackage{listings}
\numberwithin{equation}{section}

\newtheorem{assumption}{Assumption}[section]
\newtheorem{theorem}{Theorem}[section]
\newtheorem{lemma}{Lemma}[section]
\newtheorem{prop}{Proposition}[section]

\begin{document}
\title{A Runge-Kutta type scheme for nonlinear stochastic partial differential equations with multiplicative trace class noise \footnotemark[1]}
       \author{
        Xiaojie Wang \footnotemark[2] \quad  Siqing Gan  \\
       {\small School of Mathematical and Statistics,
       Central South University,}\\
      {\small Changsha 410075, Hunan,  PR China } }

       \date{}
       \maketitle

       \footnotetext{\footnotemark[1] This work was supported by NSF of China
       (No.11171352) and Hunan Provincial Innovation Foundation for Postgraduate (NO.CX2010B118). The first author would like to express his deep gratitude to Prof.P.E.Kloeden for his kind help during the author's stay in Goethe University of Frankfurt am Main. Thanks also go to Dr.A.Jentzen for very useful discussions, which were made possible through a grant from China Scholarship Council (CSC). }
        \footnotetext{\footnotemark[2]Corresponding author: x.j.wang7@gmail.com}
       \begin{abstract}
          {\rm\small In this paper a new Runge-Kutta type scheme is introduced for nonlinear stochastic partial differential equations (SPDEs) with multiplicative trace class noise. The proposed scheme converges with respect to the computational effort with a higher order than the well-known linear implicit Euler scheme. In comparison to the infinite dimensional analog of Milstein type scheme recently proposed in [Jentzen $\&$ R\"{o}ckner (2012); A Milstein scheme for SPDEs, Arxiv preprint arXiv:1001.2751v4], our scheme is easier to   implement and needs less computational effort due to avoiding the derivative of the diffusion function. The new scheme can be regarded as an infinite dimensional analog of Runge-Kutta method for finite dimensional stochastic ordinary differential equations (SODEs). Numerical examples are reported to support the theoretical results.  }\\


\end{abstract}

\section{Introduction}\label{sec:introduction}

In the last two decades, much progress has been made in developing
numerical schemes for stochastic partial differential equations (SPDEs), see, e.g.,
\cite{GK96,GI99,JK09a,JR10b,KLNS11,KLL10,LR04,LT10,ST99,ST03,Yan05},
and an extensive list of references can be found in the review article \cite{JK09}. In this article we are concerned with strong approximations (see Section
9.3 in \cite{KP92}) to nonlinear SPDEs of evolutionary type. For simplicity of presentation, we concentrate on the following
SPDE in this introductory section and refer to Section \ref{sec:Examples} for multi-dimensional space case. To be precise, we consider a parabolic SPDE with multiplicative trace class noise as
\begin{equation}\label{eq:spde}
\begin{split}
\left\{
    \begin{array}{lll} dX_t(x) = \Big[ k \frac{\partial^2}{\partial x^2}X_t(x) + f(x, X_t(x))\Big] dt + g(x, X_t(x)) d W_t(x), \quad 0 \leq t \leq T, \\
     X_t(0)=X_t(1) = 0, \\
     X_0(x) = \xi(x), \quad x \in (0,1).
    \end{array}\right.
\end{split}
\end{equation}
Here $f,g : (0,1)\times \mathbb{R} \rightarrow \mathbb{R}$ are two
appropriate smooth and regular functions with globally bounded derivatives,
and $T$ is a positive constant. Let $H=L^2((0,1),\mathbb{R})$ and let $(\Omega,\mathcal {F},\mathbb{P})$ be a probability space with a normal filtration
$\{\mathcal
{F}_t\}_{0\leq t\leq T}$. Moreover, let $W:[0,T]\times \Omega \rightarrow H$ be a standard $Q$-Wiener process with respect to $\{\mathcal
{F}_t\}_{0\leq t\leq T}$, with a trace class operator $Q:\! H \rightarrow H$.
We assume that $\eta_j, j \in \mathbb{N}$ is an orthonormal basis of $H$ consisting of eigenfunctions of $Q$ such that $Q \eta_j= \mu_j\eta_j , j\in \mathbb{N}$.


The problem \eqref{eq:spde} can be formulated in an abstract form
\begin{equation}\label{eq:abstract_spde}
dX_t = \left( AX_t + F(X_t)\right) dt + G(X_t) d W_t, \quad X_0 = \xi,
\end{equation}
where $A: \! \mathcal{D}(A) \subset H \rightarrow H$ is the Laplacian with Dirichlet boundary conditions times the constant $k >0$ and $F: \!H \rightarrow H$ and $G: \! H \rightarrow HS(U_0,H)$ are, respectively, given by $(F(v))(x) = f(x, v(x))$ and $(G(v)u)(x) = g(x,v(x))\times u(x)$ for all $x \in (0,1), v \in H$ and all $u \in U_0$. Here $U_0 = Q^{\frac{1}{2}}(H)$ and $HS(U_0,H)$ denotes the space of Hilbert-Schmidt operators from $U_0$ to $H$ (see Section \ref{sec:settings} for more details). Under the assumptions above, the SPDE \eqref{eq:abstract_spde} has a unique mild solution  with continuous sample path (see, Proposition \ref{lem:existence}), given by
\begin{equation}\label{eq:mild}
    X_t
    =
    S(t) \xi
    +
    \int_0^t S(t-s) F(X_s) \, ds
    +
    \int_0^t S(t-s) G(X_s) \, dW_s, \quad  \mathbb{P} \mbox{-a.s.},
\end{equation}
where we denote by $S(t):=e^{At}, t \geq 0$ the semigroup generated by the operator $A$.

Now we are interested in the strong approximation problem  of the SPDE \eqref{eq:spde}. More formally, we want to design a numerical approximation $Y: \Omega \rightarrow H$ such that
\begin{equation}\label{eq:strong_error}
\Big(\mathbb{E} \Big\|X_T-Y\Big\|_H^2\Big)^{\frac{1}{2}} := \Big(\mathbb{E} \Big[\int_{(0,1)}|X_T(x)-Y(x)|^2dx \Big]\Big)^{\frac{1}{2}} < \varepsilon
\end{equation}
holds for a given precision $\varepsilon > 0$ with the least possible computational effort. To implement the numerical approximation on a
computer, one has to discretize both the time interval $[0,T]$ and the infinite dimensional space $H$. In this article we consider spectral Galerkin method for spatial discretization and difference method for temporal discretization. A simple fully discretization for \eqref{eq:spde} is the linear implicit Euler scheme combined with spectral Galerkin method given by
$ \bar{Y}_0^{N,M,K} = P_N( \xi ) $
and for $ m =0,1,\ldots,M-1$
\begin{align}
\label{eq:Escheme}
&
  \bar{Y}_{m+1}^{N,M,K}
  =
  P_N \left(  \Big(I -h A\Big)^{-1}
  \Big(
    \bar{Y}_m^{N,M,K}
    +
    h \,f\!\left(\cdot, \bar{Y}_m^{N,M,K} \right)
    +
    g\!\left(\cdot, \bar{Y}_m^{N,M,K}
    \right) \!\times
    \Delta W_m^{M,K}
  \Big)\right),
\end{align}
where we used the notations $\varphi(\cdot, v): (0,1) \rightarrow \mathbb{R}$ and $v \times w: (0,1) \rightarrow \mathbb{R}$ given by
\begin{equation}\label{eq:notations1}
(\varphi(\cdot, v))(x) = \varphi(x, v(x)),  \quad (v \times w)(x) = v(x)\times w(x)
\end{equation}
for all $x \in (0,1)$ and all functions $v, w: (0,1) \rightarrow \mathbb{R}$ and $\varphi : (0,1) \times \mathbb{R} \rightarrow \mathbb{R}$.
Here and throughout this article,
$ h = \frac{T} {M}, M \in \mathbb{N} $ is the time stepsize and the increment $\Delta W^{M,K}_m(\omega) :=
W^K_{ \frac{ (m + 1) T }{ M } }(\omega) -
W^K_{ \frac{ m T }{ M } }(\omega)$ is given by \eqref{Wiener_Incr}.
The linear projection operator $P_N: H \rightarrow H$ is defined by \eqref{def:PN}.
Here $\{e_i\}_{i \in \mathbb{N}}$ is an orthonormal basis of $U$ consisting of eigenfunctions of $A$. It is worthwhile to point out that the scheme \eqref{eq:Escheme} for \eqref{eq:spde} is easy to implement (see Figure 2 in \cite{JR10b} for the matlab code).

Recently, Jentzen and R\"{o}ckner \cite{JR10b} introduced an infinite dimensional analog of Milstein type scheme  for \eqref{eq:spde},
given by
$ \tilde{Y}_0^{N,M,K} = P_N( \xi ) $
and for $m=0,1,...,M-1$
\begin{align}
\label{eq:CMscheme}
&
  \tilde{Y}_{m+1}^{N,M,K}
  =
  P_N \left( S (h)
  \Bigg(
    \tilde{Y}_m^{N,M,K}
    +
    h \,
    f\!\left(\cdot, \tilde{Y}_m^{N,M,K} \right)
    +
    g\!\left(\cdot, \tilde{Y}_m^{N,M,K}
    \right) \!\times
    \Delta W_m^{M,K}\right.
\\&\left.
\nonumber
    +
    \frac{1}{2}
    \Big(\frac{\partial }{\partial y}g\Big)\left( \cdot,\tilde{Y}_m^{N,M,K} \right) \times
    g\left( \cdot, \tilde{Y}_m^{N,M,K} \right) \times \Big(
      (\Delta W_m^{M,K})^2
    - h
      \sum_{j=1}^K
        \mu_j (\eta_j)^2 \Big)\Bigg)\right).
\end{align}
Here, apart from \eqref{eq:notations1}, we also used the notations $v^2: (0,1) \rightarrow \mathbb{R}$ and $(\frac{\partial }{\partial y}g)(\cdot, v): (0,1) \rightarrow \mathbb{R}$ given by
\begin{equation}\nonumber
(v^2)(x) = (v(x))^2, \quad \left(\frac{\partial }{\partial y}g\right)(\cdot, v)(x) = \frac{\partial }{\partial y}g(x, v(x))
\end{equation}
for all $x \in (0,1)$ and all functions $v,w: (0,1) \rightarrow \mathbb{R}$. For continuously differentiable function $g : (0,1) \times \mathbb{R} \rightarrow \mathbb{R}$, $\frac{\partial }{\partial y}g: (0,1) \times \mathbb{R} \rightarrow \mathbb{R} $ is a partial derivative of $g$ with respect to the second variable.  On the one hand, it is also easy to implement the scheme \eqref{eq:CMscheme} for the SPDE \eqref{eq:spde} (see Figure 3 in \cite{JR10b} for the matlab code). On the other hand, the scheme \eqref{eq:CMscheme} gives a break of complexity of the numerical approximation of nonlinear SPDEs with multiplicative trace class noise.  For example, the scheme \eqref{eq:Escheme} for the first test example in Section \ref{sec:Numer-Experi} can only achieve overall convergence order $\frac{3}{8}-$ while the scheme \eqref{eq:CMscheme} possesses overall convergence order  $\frac{1}{2}-$ (here and below we write
$b-$ for the convergence order if the convergence order is higher than $b-\epsilon$ for every arbitrarily small $0< \epsilon <b$).

Note that the scheme \eqref{eq:CMscheme} can also be adapted to solve a
SPDE system \eqref{eq:spde}
with $X_t(x), \xi(x) \in \mathbb{R}^n$, and $f,g : (0,1)\times \mathbb{R}^n \rightarrow \mathbb{R}^n$ for $n \in \mathbb{N}$. In this situation, $\frac{\partial g}{\partial y}$ in \eqref{eq:CMscheme} is interpreted as the Jacobian matrix of $g$. Therefore, to implement \eqref{eq:CMscheme} one needs to calculate the Jacobian exactly, which may be difficult, and to evaluate it at each time step, which may be expensive. To save computational cost in this sense, we will take a Runge-Kutta type scheme to avoid computing the Jacobian. Following this idea, in this article we aim at constructing a high strong order Runge-Kutta method for \eqref{eq:spde}. For simplicity, we only consider a scalar SPDE. But our work can be easily extended to a SPDE system with scalar noise as described above.

One approach for deriving a Runge-Kutta method is to replace the partial derivative in the approximations \eqref{eq:CMscheme} by difference, and this leads to a derivative-free scheme, given by $Y_0^{N,M,K} = P_N( \xi )$ and for $m=0,1,...,M-1$
\begin{align}\label{eq:concrete_scheme}
&Y_{m+1}^{N,M,K} = P_N \Bigg( S(h) \Bigg( Y_m^{N,M,K} + h \, f(\cdot, Y_m^{N,M,K}) + g\left(\cdot, Y_m^{N,M,K}\right)\times \Delta W_m^{M,K}
\\  \nonumber & +
\frac{1}{2\sqrt{h}} \left[ g\Big(\cdot, Y_m^{N,M,K} + \sqrt{h}\,g(\cdot, Y_m^{N,M,K})\Big)-g(\cdot, Y_m^{N,M,K}) \right]\times \Big((\Delta W^{M,K}_m)^2-h\sum_{j=1}^{K} \mu_j(\eta_j)^2\Big) \Bigg) \Bigg),
\end{align}
where the function $g\Big(\cdot, v+\sqrt{h}\,g(\cdot, v)\Big)-g(\cdot, v): (0,1) \rightarrow \mathbb{R}$ is defined by
\begin{equation}\nonumber
\left[g\Big(\cdot, v+\sqrt{h}\,g(\cdot, v)\Big)-g\Big(\cdot, v\Big)\right](x) = g\Big(x, v(x)+\sqrt{h}\,g(x, v(x))\Big)-g\Big(x, v(x)\Big)
\end{equation}
for all functions $v: (0,1) \rightarrow \mathbb{R}, g : (0,1) \times \mathbb{R} \rightarrow \mathbb{R}$.
A natural question thus arises as to whether such replacement maintains the high
convergence order of \eqref{eq:CMscheme}. In this paper we give a positive answer
and prove that the new scheme not only maintains the high convergence order,
but also reduces  computational cost. Similarly to the scheme \eqref{eq:CMscheme}, the numerical method \eqref{eq:concrete_scheme} can be simulated quite easily (see Figure \ref{code} in Section \ref{sec:Numer-Experi} for the implementation code).

Now we take a closer look at schemes \eqref{eq:CMscheme} and \eqref{eq:concrete_scheme}. For each step, the Milstein type scheme \eqref{eq:CMscheme} requires one evaluation of $f$, one evaluation of $g$ and one evaluation of the partial derivative $\frac{\partial }{\partial y}g$. In contrast,  the new Runge-Kutta type scheme \eqref{eq:concrete_scheme} needs one evaluation of $f$ and two evaluations of $g$, but no evaluation of the partial derivative $\frac{\partial }{\partial y}g$ at each step. Thus the Runge-Kutta type scheme \eqref{eq:concrete_scheme} is easier to implement than the Milstein type scheme \eqref{eq:CMscheme}. The main result (Theorem \ref{thm:mainresult}) shows that the new scheme \eqref{eq:concrete_scheme} maintains the high convergence order of scheme \eqref{eq:CMscheme}. Numerical results in Section \ref{sec:Numer-Experi} demonstrate that the schemes \eqref{eq:CMscheme} and \eqref{eq:concrete_scheme} produce nearly the same approximation errors. Even if we neglect the effort for the calculation of the partial derivative, the runtime for one path simulation of the Runge-Kutta scheme \eqref{eq:concrete_scheme} applied to \eqref{eq:spde} is less than that for simulating the scheme \eqref{eq:CMscheme}. Take the first test problem in Section \ref{sec:Numer-Experi} for example, in the case of $N=256$, one path simulation of the Runge-Kutta scheme \eqref{eq:concrete_scheme} costs 60.312000 seconds, while it needs 77.016000 seconds to simulate the scheme \eqref{eq:CMscheme} (see Table \ref{table1}). This occurs due to the fact that evaluation of the partial derivative $(\frac{\partial }{\partial y}g)(x,y)$ costs more time than evaluation of the function $g(x,y)$.  Summarizing, the Runge-Kutta type scheme \eqref{eq:concrete_scheme} is easier to implement and needs  less computational effort than the Milstein type scheme \eqref{eq:CMscheme}. Moreover, the derivative-free scheme \eqref{eq:concrete_scheme} can be regarded as an infinite dimensional analog of the Runge-Kutta method (6) in \cite{BB96} for finite SODEs.  To the best of our knowledge, this is the very first paper to introduce a Runge-Kutta method for nonlinear SPDEs with multiplicative noise. We also mention that constructing higher order Runge-Kutta methods and developing more systematic way to derive Runge-Kutta methods for SPDEs are of great interest and will be our future work.

The rest of this paper is organized as follows. In the next section, we put everything into an abstract framework and state the main convergence result of this article. In Section \ref{sec:Examples} we give examples fulfilling the assumptions in the previous section. A detailed proof of the main convergence result is elaborated in Section \ref{sec:proof}. Finally, we illustrate how to implement the proposed scheme and  present some numerical examples to support our theoretical results.

\section{Abstract framework and main result}\label{sec:settings}
In this section we focus on the abstract framework \eqref{eq:abstract_spde}
and adopt the following setting and assumptions.

Let
$
  \left(
    H,
    \left< \cdot , \cdot \right>_H,
    \left\| \cdot \right\|_H
  \right) $
and
$
  \left(
    U,
    \left< \cdot , \cdot \right>_U,
    \left\| \cdot \right\|_U
  \right)
$
be two separable Hilbert spaces. By $L(U,H)$ and
$L^{(2)}(U, H)$ we denote the space of linear bounded operator
from $U$ to $H$ and from $U \times U$ to $H$, respectively.
For short, we write $L(H)$ instead of $L(H,H)$.  We also introduce a space of Hilbert-Schmidt operators. An operator $\Gamma \in L(U,H)$ belongs to the Hilbert-Schmidt operator space $HS(U,H)$, if for any orthonormal basis $\{\psi_k\}_{k=1}^{\infty}$ of $U$ the sum
\begin{equation}\nonumber
\|\Gamma\|_{HS(U,H)}^2 := \sum_{k=1}^{\infty} \|\Gamma \psi_k\|_{H}^2
\end{equation}
is finite and is independent of the choice of the orthonormal basis. The quantity $\|\Gamma\|_{HS(U,H)}$ is called the Hilbert-Schmidt norm of $\Gamma$. Similarly, we can define an Hilbert-Schmidt operator space $HS^{(2)}(U,H)$
from $U \times U$ to $H$. We refer to Chow \cite{Cpl07}, Da Prato and Zabczyk \cite{DZ92}, Pr\'{e}v\^{o}t and R\"{o}ckner \cite{PR07} for
details on these spaces and their properties.

Moreover, let $Q \in L(U) $ be symmetric, nonnegative and with finite trace, i.e.,
\begin{equation} \label{eq:trQ}
Tr(Q) < \infty.
\end{equation}
Suppose that $ \mathcal{J} $ is a finite or countable set, and let
$ \left( \eta_j \right)_{ j \in \mathcal{J} }
\subset U $ be an orthonormal basis of $U$ consisting of eigenfunctions of $Q: U \rightarrow U$ such that $Q \eta_j = \mu_j \eta_j , j\in \mathcal{J}$. Denote by
$
  \left(
    U_0,
    \left< \cdot , \cdot \right>_{ U_0 },
    \left\| \cdot \right\|_{ U_0 }
  \right)
$
the separable Hilbert space
$ U_0 := Q^{ \frac{1}{2} }( U ) $
with $ \left< v, w \right>_{ U_0 } =
\big< Q^{ -\frac{1}{2} } v, Q^{ -\frac{1}{2} } w \big>_U $
for all $ v, w \in U_0 $, where $Q^{ -\frac{1}{2} }$ is the pseudo inverse of $Q^{\frac{1}{2}}$ (see, e.g., Section 2.3.2 in \cite{PR07}). We can obtain that
\begin{equation}\label{eq:HS}
\|\Gamma\|_{HS(U_0, H)} = \left\|\Gamma \circ Q^{ \frac{1}{2} } \right\|_{HS(U, H)} \quad \mbox{for} \quad \Gamma \in HS(U_0, H).
\end{equation}

Further, we assume that $(\Omega,\mathcal {F},\mathbb{P})$ is a probability space with a normal filtration
$\{\mathcal
{F}_t\}_{0\leq t\leq T}$ and $W:[0,T]\times \Omega \rightarrow U$ is a standard $Q$-Wiener process with respect to $\{\mathcal
{F}_t\}_{0\leq t\leq T}$ and has the representation \cite[Proposition 2.1.10]{PR07}
\begin{equation}\label{Wiener_Process_Rep}
W_t( \omega ) :=
\sum_{
  \substack{
    j \in \mathcal{J}
  \\
    \mu_j \neq 0
  }
}
\sqrt{\mu_j} \beta_t^j(\omega) \eta_j,
\end{equation}
where $(\beta^j_t)_{j \in \mathcal{J}, \mu_j \neq 0}$ for $t \in [0,T]$ are independent real-valued Brownian motions on the probability space $\left( \Omega, \mathcal{F}, \{\mathcal
{F}_t\}_{0\leq t\leq T}, \mathbb{P} \right)$. Now we make the following assumptions.

\begin{assumption}[Linear operator $A$]\label{semigroup}
Let $ \mathcal{I} $ be a
finite or countable set
and let $ \left( \lambda_i
\right)_{i \in \mathcal{I} } $
be a family of real numbers with
$ \inf_{ i \in \mathcal{I} } \lambda_i \in (0,\infty) $.
Further, let
$ (e_i)_{ i \in \mathcal{I} } $
be an orthonormal basis of $H$ and
let $ -A : \mathcal{D}(-A) \subset H \rightarrow H $
be a linear operator such that
\begin{equation}\label{A_Sper}
  -Av
  =
  \sum_{ i \in \mathcal{I} }
    \lambda_i
    \left<
      e_i, v
    \right>_H
    e_i
\end{equation}
for all
$v \in \mathcal{D}(-A)
  :=
  \left\{
    v \in H
    \big|
    \sum_{ i \in \mathcal{I} }
      \left|
        \lambda_i
      \right|^2
      \left|
        \left<
          e_i, v
        \right>_H
      \right|^2
    <
    \infty
  \right\}
$.
\end{assumption}
Here and below we denote by $S(t):=e^{At}, t \geq 0$ the semigroup generated by the operator $A$. By $ V_r := \mathcal{D}\left( \left( - A \right)^r \right), \: r \geq 0 $
equipped with the norm
$ \left\| v \right\|_{ V_r }
:= \left\| \left( - A \right)^r v \right\|_H $
we denote
the $ \mathbb{R} $-Hilbert spaces
of domains of fractional powers
of the linear operator $ -A$.

\begin{assumption}[Drift
coefficient $F$]\label{drift}
For $\beta \in [0,1) $, we assume that $ F : V_{\beta} \rightarrow H $ is a twice continuously Fr\'{e}chet differentiable mapping with
$
  \sup_{ v \in V_{\beta} }
  \left\| F'(v) \right\|_{ L(H) }
  < \infty
$
and with
$
  \sup_{ v \in V_{\beta} }
  \left\| F''(v) \right\|_{ L^{(2)}(V_{\beta},H) }
  < \infty
$.
\end{assumption}

\begin{assumption}[Diffusion
coefficient $G$]\label{diffusion}
Let $ G : V_{ \beta } \rightarrow HS(U_0,H) $
be a twice continuously Fr\'{e}chet differentiable
mapping such that $ \sup_{ v \in V_{\beta} } \!
\| G'( v )
\|_{ L( H, HS( U_0, H ) ) } \!
< \infty $ and\\
$ \sup_{ v \in V_{\beta} } \!
\| G''( v )
\|_{ L^{(2)}( V_{\beta}, HS( U_0, H ) ) }
$
$<$
$\infty $.
Moreover,
let
$ \alpha, c \in (0,\infty) $,
$ \delta, \vartheta
\in ( 0, \frac{1}{2} )$ with $ \beta \leq \delta + \frac{1}{2} $, and
$ \gamma \in [ \max(\delta,\beta),
\delta + \frac{1}{2} ) $.
Suppose that $ G( V_{ \delta } ) $
$\subset$ $ HS(U_0, V_{ \delta } ) $
and
\begin{equation}\label{condition_AA}
  \left\| G(u) \right\|_{ HS( U_0, V_{ \delta } ) }
  \leq
  c \left( 1 + \left\| u \right\|_{ V_{ \delta } } \right) ,
\end{equation}
\begin{equation}\label{condition_A}
    \left\|
      G'\!\left( v \right) G\!\left( v \right)
      -
      G'\!\left( w \right) G\!\left( w \right)
    \right\|_{ HS^{ (2) }( U_0, H ) }
  \leq c
    \left\| v - w \right\|_H ,
\end{equation}
\begin{equation}\label{condition_AAA}
  \left\|
    \left( - A \right)^{ - \vartheta } \!
    G(v) Q^{ - \alpha }
  \right\|_{ HS( U_0, H ) }
  \leq
  c \left( 1 + \left\| v \right\|_{ V_{ \gamma } } \right)
\end{equation}
hold for all $ u \in V_{ \delta } $
and $ v, w \in V_{ \gamma } $.
Furthermore, let the bilinear Hilbert-Schmidt
operator $ G'(v) G(v) \in HS^{(2)}(U_0, H) $
be symmetric for all $ v \in V_{ \beta } $.
\end{assumption}

Note that
the operator $ G'(v) G(v) : U_0 \times U_0 \rightarrow H $
given by
\begin{equation}
\label{def:bilinear}
  \Big( G'(v) G(v) \Big)
  ( u, \tilde{u} )
  =
  \Big(
    G'(v) \big( G(v) u \big)
  \Big)
  ( \tilde{u} )
\end{equation}
for all $ u, \tilde{u} \in U_0 $ is
a bilinear Hilbert-Schmidt operator in $ HS^{(2)}( U_0, H ) $
for all $ v \in V_{\beta} $.
The assumed symmetry of
$ G'(v)G(v) \in HS^{ (2) }( U_0, H ) $
thus reads as \cite[Remark 1]{JR10b}
\begin{equation}
\label{eq:commutativity}
  \Big(
    G'(v) \big( G(v) u \big)
  \Big)
  ( \tilde{u} )
=
  \Big(
    G'(v) \big( G(v) \tilde{u} \big)
  \Big)
  ( u )
\end{equation}
for all $ u, \tilde{u} \in U_0 $
and all $ v \in V_{ \beta } $.
We also mention that
\eqref{eq:commutativity} is the
abstract (infinite dimensional) analog
of the commutativity condition (10.3.13) in \cite{KP92}.
Although the commutativity condition (10.3.13) in \cite{KP92}
is seldom fulfilled for
finite dimensional SODEs,
\eqref{eq:commutativity}
is naturally met for SPDE \eqref{eq:spde}.

\begin{assumption}[Initial value $\xi$]\label{initial}
Let $ \xi : \Omega \rightarrow V_{\gamma} $ be
an $ \mathcal{F}_0
$/$ \mathcal{B}\left(
V_{\gamma} \right) $-measurable mapping
with $ \mathbb{E} \left\| \xi \right\|^4_{ V_{\gamma} }
< \infty $.
\end{assumption}
We remark that a mapping $ \xi : \Omega \rightarrow V_{\gamma} $ is called $\mathcal{F}_0/\mathcal{B}(V_{\gamma})$ measurable if it is a measurable mapping from the measurable space $(\Omega, \mathcal{F}_0)$ to the measurable space $(V_{\gamma}, \mathcal{B}(V_{\gamma}))$. Here $\mathcal{B}(V_{\gamma})$ denotes the Borel $\sigma$-field of $V_{\gamma}$. The assumptions above are sufficient to guarantee the existence of a
unique mild solution of the SPDE~\eqref{eq:abstract_spde} \cite[Theorem~1]{JR10a}.
\begin{prop}[Existence of the mild solution]\label{lem:existence}
Let Assumptions~\ref{semigroup}-\ref{initial} and condition \eqref{eq:trQ} be fulfilled.
Then there exists an up to modifications
unique predictable stochastic process
$X : [0,T] \times \Omega \rightarrow V_{\gamma} $,
which fulfills
$
  \sup\limits_{ t \in [0,T] }
  \mathbb{E}
    \left\|
      X_t
    \right\|_{ V_{ \gamma } }^4
  < \infty
$,
$
  \sup\limits_{ t \in [0,T] }
  \mathbb{E}
    \left\|
      G(X_t)
    \right\|_{ HS(U_0,V_{ \delta }) }^4
  < \infty
$ and
\begin{equation}\label{eq:SPDE}
    X_t
    =
    S(t) \xi
    +
    \int_0^t S(t-s) F(X_s) \, ds
    +
    \int_0^t S(t-s) G(X_s) \, dW_s, \quad \mathbb{P}\mbox{-a.s.}
\end{equation}
for all $t \in [0,T]$.
\end{prop}
Let $ \left( \mathcal{I}_N \right)_{ N \in \mathbb{N} } $
and $ \left( \mathcal{J}_K \right)_{ K \in \mathbb{N} } $ be
sequences of finite subsets of $ \mathcal{I} $ and $ \mathcal{J} $,
respectively.
For $ N \in \mathbb{N} $ we define the linear projection operators
$ P_N : H \rightarrow H $, by
\begin{equation}\label{def:PN}
P_N( v ) := \sum_{ i \in \mathcal{I}_N }
\left< e_i, v \right>_H e_i, \quad v \in H.
\end{equation}
Furthermore, for all $ K \in \mathbb{N} $ we define Wiener processes $ W^K : [0,T] \times \Omega \rightarrow U_0 $ by
\begin{equation}\label{Wiener_Process}
W_t^K( \omega ) :=
\sum_{
  \substack{
    j \in \mathcal{J}_K
  \\
    \mu_j \neq 0
  }
}
\sqrt{\mu_j} \beta_t^j(\omega) \eta_j, \quad t \in [0,T], \: \: \omega \in \Omega.
\end{equation}
Here $(\beta^j)_{j \in \mathcal{J}, \mu_j \neq 0}$ are independent real-valued Brownian motions on the probability space $\left( \Omega, \mathcal{F}, \{\mathcal
{F}_t\}_{0\leq t\leq T}, \mathbb{P} \right)$. We also use
the notations $\Delta \beta_m^j(\omega) := \beta_{t_{m+1}}^j(\omega) - \beta_{t_m}^j(\omega)$ and
\begin{equation}\label{Wiener_Incr}
\Delta W^{M,K}_m(\omega) :=
W^K_{ \frac{ (m + 1) T }{ M } }(\omega) -
W^K_{ \frac{ m T }{ M } }(\omega) = \sum_{
  \substack{
    j \in \mathcal{J}_K
  \\
    \mu_j \neq 0
  }
}
\sqrt{\mu_j} \Delta \beta_m^j(\omega) \eta_j
\end{equation}
for all
$\omega \in \Omega $ and $ m = 0, 1, \dots, M - 1$.

Subsequently, we formulate the schemes in the introduction part  in abstract form. Then the scheme \eqref{eq:CMscheme} for the problem \eqref{eq:spde} can be formulated in an abstract scheme for the abstract problem \eqref{eq:abstract_spde}, given by
$ \tilde{Y}_0^{N,M,K} = P_N( \xi ) $
and for $ m =0,1,\ldots,M-1 $
\begin{align}
\label{eq:Mscheme}
&
  \tilde{Y}_{m+1}^{N,M,K}
  =
  P_N \Bigg( S( h )
  \Bigg(
    \tilde{Y}_m^{N,M,K}
    +
    h \,
    F\!\left( \tilde{Y}_m^{N,M,K} \right)
    +
    G\!\left( \tilde{Y}_m^{N,M,K}
    \right) \!
    \Delta W_m^{M,K}
\\&
\nonumber
    +
    \frac{1}{2}
    G'\!\left( \tilde{Y}_m^{N,M,K} \right) \!
    \Big( \!
      G\!\left( \tilde{Y}_m^{N,M,K} \right) \!
      \Delta W_m^{M,K}
    \Big)
    \Delta W_m^{M,K}
    - \frac{h}{2} \!\!
      \sum_{
        \substack{
          j \in \mathcal{J}_K
        \\
          \mu_j \neq 0
        }
      }
        \!\!
        \mu_j
        G'\!\left( \tilde{Y}_m^{N,M,K} \right) \!
        \Big( \!
          G\!\left( \tilde{Y}_m^{N,M,K} \right)
          \eta_j
        \Big)
        \eta_j
  \Bigg)\Bigg).
\end{align}

Assume that Assumptions \ref{semigroup}-\ref{initial} are all satisfied, Jentzen and R\"{o}ckner \cite{JR10b} have established the strong convergence result for
the scheme \eqref{eq:Mscheme} (see Theorem 1 in \cite{JR10b}). In later development, we show that this convergence result also holds for the Runge-Kutta type method proposed here under some additional conditions. First of all, we formulate the new numerical scheme \eqref{eq:concrete_scheme} in an abstract form, given by $Y_0^{N,M,K} = P_N( \xi )$ and
for $ m =0,1,\ldots,M-1 $
\begin{align}\label{eq:abstract_scheme}
Y_{m+1}^{N,M,K} =& P_N \Bigg(S(h) \Bigg( Y_m^{N,M,K} + h \,F(Y_m^{N,M,K}) + G\left( Y_m^{N,M,K}\right)\Delta W_m^{M,K}
\\ & +\frac{1}{2} GG(Y_m^{N,M,K},h)\left(\Delta W^{M,K}_m, \Delta W^{M,K}_m \right)-\frac{h}{2}
\sum_{\begin{subarray}{ll} j\in \mathcal{J}_K \nonumber \\ \mu_j \neq 0 \end{subarray}}
\mu_j GG(Y_m^{N,M,K},h)\left(\eta_j, \eta_j\right)\Bigg)\Bigg),
\end{align}
where $ GG(v,h) : U_0 \times U_0 \rightarrow H $ is a derivative-free bilinear operator to approximate the bilinear operator $ G'(v) G(v) : U_0 \times U_0 \rightarrow H $ in \eqref{eq:Mscheme}. We define the remainder bilinear operators $GG_1(v,h): \! U_0 \times U_0 \rightarrow H $ by
\begin{equation}\label{def:BB1}
\Big( GG_1(v,h) \Big)
  ( u, \tilde{u} ) := \Big( GG(v,h) \Big)
  ( u, \tilde{u} ) -\Big(G'(v)(G(v)u)\Big)(\tilde{u}).
\end{equation}

Note that Assumptions \ref{semigroup}-\ref{initial} come from \cite{JR10b}, which are used to ensure the high strong convergence order of the Milstein-type scheme \eqref{eq:CMscheme}. In the next section, concrete conditions are given for concrete parabolic SPDE to promise these assumptions. To guarantee the high convergence order of the Runge-Kutta scheme, we shall impose additional conditions on the
two bilinear operators $ GG(v,h)$ and $GG_1(v,h)$ as follows.

\begin{assumption}[Approximation operator]\label{BB_assum}
Suppose that for any $v,w \in V_{\beta}$ there exists a constant $C_0$ independent of $h$ such that
\begin{align}\label{eq:BBLip}
\|GG(v,h)-GG(w,h)\|_{HS^{(2)}(U_0,H)}^2 \leq& \frac{C_0}{h} \|v-w\|_H^2,
\\  \|GG_1(v,h)\|_{HS^{(2)}(U_0,H)}^2 \leq & C_0 h \big(1+ \|v\|_{V_{\beta}}^4\big).
\label{eq:BB1GC}
\end{align}
\end{assumption}
In the next section we will validate the imposed conditions \eqref{eq:BBLip}-\eqref{eq:BB1GC} for parabolic SPDEs (see Proposition \ref{lem:BB1}). Armed with these assumptions, we are now to give the main result of this article.

\begin{theorem}[Main
result]\label{thm:mainresult}
Suppose that Assumptions~\ref{semigroup}-\ref{BB_assum} and \eqref{eq:trQ} are fulfilled. Then there is
a constant $ C >0 $ independent of $h$ such that
\begin{multline}\label{eq:mainresult}
    \sup_{0\leq m\leq M }\left(
      \mathbb{E}
      \left\|
        X_{ mh }
        - Y^{N,M,K}_m
      \right\|_H^2
    \right)^\frac{1}{2}
    \leq
    C
  \bigg(
      \Big(
        \inf_{
          i \in \mathcal{I} \backslash
          \mathcal{I}_N
        }
        \lambda_i
      \Big)^{ - \gamma }
    +
     \Big(
      \sup_{ j \in \mathcal{J} \backslash \mathcal{J}_K }
      \mu_j
      \Big)^{ \alpha }
    +
      M^{ - \min\left( 2 \left( \gamma - \beta \right) ,
      \gamma \right) }
  \bigg).
\end{multline}
\end{theorem}
The detailed proof is postponed to Section \ref{sec:proof}. Theorem \ref{thm:mainresult} indicates that the approximation error in \eqref{eq:mainresult} is composed of three parts. The first term
$\left(
        \inf_{
          i \in \mathcal{I} \backslash
          \mathcal{I}_N
        }
        \lambda_i
      \right)^{ - \gamma }$
arises due to spatial discretization. The second term $\big(
      \sup_{ j \in \mathcal{J} \backslash \mathcal{J}_K }
      \mu_j
      \big)^{ \alpha }$ comes from truncation of the expansion of the noise $W_t$. The third term $M^{ - \min\left( 2 \left( \gamma - \beta \right) ,
      \gamma \right) }$ corresponds to the temporal discretization error.

\section{Parabolic SPDEs}\label{sec:Examples}


In this section we give concrete parabolic SPDE examples falling into the abstract framework in Section \ref{sec:settings}. Let $d \in \{1,2,3\}$ and let $H=U=L^2((0,1)^d, \mathbb{R})$ be the Hilbert space with the scalar product and the norm, respectively, given by
$$
\langle v,w \rangle_H = \int_{(0,1)^d} v(x)\times w(x) dx \quad\quad \mbox{and} \quad\quad \|v\|_H = \left(\int_{(0,1)^d} |v(x)|^2 dx\right)^{\frac{1}{2}}.
$$
For the continuous function $v: (0,1)^d \rightarrow \mathbb{R}$,  we define two norms
$$
\|v\|_{C((0,1)^d;\mathbb{R})}:= \sup_{x\in (0,1)^d} |v(x)|
$$
and
$$
\|v\|_{C^r((0,1)^d;\mathbb{R})}:= \sup_{x\in (0,1)^d} |v(x)|
\:+ \sup_{x,y\in (0,1)^d, x \neq y} \frac{|v(x)-v(y)|}{\|x-y\|_{\mathbb{R}^d}^r},
$$
where the Euclidean norm $\|x\|_{\mathbb{R}^d} := (|x_1|^2+...+|x_d|^2 )^{\frac{1}{2}}$ was used for $x = (x_1,...,x_d) \in \mathbb{R}^d$.
First of all, we assume that for some constants $0<\rho<1, c >0$, the eigenfunctions $\eta_j, j \in  \mathcal{J} $ of the covariance operator $Q$ are continuous and satisfy
\begin{equation}\label{eq:gj}
\sup_{ j\in \mathcal{J}} \|\eta_j\|_{C((0,1)^d;\mathbb{R})} \leq c, \quad \sum_{j\in \mathcal{J}} \mu_j \|\eta_j\|^2_{C^{\rho}((0,1)^d;\mathbb{R})} \leq c.
\end{equation}
Here and below $c$ is a generic constant, which may be different in different places.

For the linear operator $A$ in Assumption \ref{semigroup}, let $\mathcal{I}=\mathbb{N}^d$ and $A:= k \Delta$ with $k > 0$ be the Laplacian times a constant with Dirichlet boundary condition, i.e.,
$
A v = k \Delta v = k \left( \sum_{j=1}^{d} \frac{\partial^2}{\partial x_j^2} \right)v
$
for $v \in D(-A)$. Then \eqref{A_Sper} in Assumption \ref{semigroup} holds with
\begin{equation}\nonumber
e_i(x)= 2^{\frac{d}{2}} \prod_{j=1}^d\sin (i_j \pi x_j), \quad \lambda_i = k \pi^2 \sum_{j=1}^d (i_j)^2
\end{equation}
for all $x = (x_1,...,x_d) \in (0,1)^d$ and all $i=(i_1,...,i_d)\in \mathbb{N}^d$. Here we set $\mathcal{I}_N = \{1,...,N\}^d$.

For the drift coefficient $F$ in Assumption \ref{drift}, set
$$
\beta = \frac{d}{5}, \quad \mbox{for} \quad d =1,2,3,
$$
and  let $f: (0,1)^d \times \mathbb{R} \rightarrow \mathbb{R}$ be a twice continuously differentiable function such that
\begin{equation}\label{eq:drift}
\begin{split}
\int_{(0,1)^d}|f(x,0)|^2 dx \leq& c, \quad \Big|\Big(\frac{\partial^n }{\partial y^n}f\Big) (x,y) \Big| \leq c, \quad n=1,2
\end{split}
\end{equation}
for $ x \in (0,1)^d, y \in \mathbb{R}$. Then we define the operator $F: V_{\beta} \rightarrow H$ as
\begin{equation}\label{def:F}
(F(v))(x) := f(x, v(x))
\end{equation}
for all $x \in(0,1)^d$ and all $v \in V_{\beta}$. 

For the diffusion coefficient $G$ in Assumption \ref{diffusion}, let $g: (0,1)^d \times \mathbb{R} \rightarrow \mathbb{R}$ be a twice continuously differentiable function with
\begin{equation}\label{eq:diffusion}
\begin{split}
|g(x,0)| \leq c, \quad \Big|\Big(\frac{\partial^n }{\partial y^n}g\Big) (x,y) \Big| \leq c, \quad \Big\|\Big(\frac{\partial }{\partial x}g\Big)(x,y)\Big\|_{L(\mathbb{R}^d; \mathbb{R})} \leq c, \: n=1,2,
\end{split}
\end{equation}
and
\begin{equation}\label{eq:diffusion2}
\left|\Big(\frac{\partial }{\partial y}g\Big)(x,y)\, g(x,y)- \Big(\frac{\partial }{\partial y}g\Big)(x,z)\, g(x,z) \right| \leq c |y-z|
\end{equation}
for $x \in (0,1)^d, y,z \in \mathbb{R}$. Here $\|\cdot\|_{L(\mathbb{R}^d; \mathbb{R})}$ is the usual operator norm. Then let the operator $G: V_{\beta} \rightarrow HS(U_0,H)$ be given by
$$
(G(v)u)(x) := g(x, v(x)) \times u(x)
$$
for all $x \in(0,1)^d, v \in V_{\beta}$ and $u \in U_0 \subset U=H$. Therefore the bilinear operator $G'(v)G(v)$ in the abstract scheme \eqref{eq:Mscheme} is here given by
\begin{equation}
\label{def:con_bilinear}
  \Big( G'(v) G(v) \Big)
  ( u, \tilde{u} )(x)
  = \Big(\frac{\partial }{\partial y}g\Big)(x,v(x))\, g(x,v(x))\times u(x) \times \tilde{u}(x).
\end{equation}

In regard to the initial value in Assumption \ref{initial}, let $x_0: [0,1]^d \rightarrow \mathbb{R}$ be a twice continuously differentiable function with $x_0|_{\partial (0,1)^d} \equiv 0$. Then let the initial value be given by $\xi(\omega)=x_0$ for all $\omega \in \Omega$.

It is shown in \cite[Section 4]{JR10b} that the linear operator $A$, the drift coefficient $F$, the diffusion coefficient $G$ and the initial value $\xi$ defined as above satisfy all conditions in Assumption \ref{semigroup}-\ref{initial}, except for \eqref{condition_AA} and \eqref{condition_AAA} in Assumption \ref{diffusion}, which will be verified for some concrete examples in Section \ref{sec:Numer-Experi}.

In the setting above, SPDE \eqref{eq:abstract_spde} reduces to a parabolic SPDE as
\begin{equation}\label{eq:example_spde}
\begin{split}
dX_t(x) = \left[ k \left(\sum_{j=1}^d \frac{\partial^2}{\partial x_j^2}\right) X_t(x) + f(x, X_t(x))\right] dt + g(x, X_t(x)) d W_t(x)
\end{split}
\end{equation}
with $X_t|_{\partial (0,1)^d} \equiv 0$ and $X_0(x) = x_0(x)$ for $x \in (0,1)^d$ and $t \in [0,T]$. For \eqref{eq:example_spde}, linear implicit Euler method and Milstein type method  take the same form as \eqref{eq:Escheme} and \eqref{eq:CMscheme}, respectively. If we introduce a bilinear operator $ GG(v,h) : U_0 \times U_0 \rightarrow H $ approximating the bilinear operator $G'(v)G(v)$ in the scheme \eqref{eq:abstract_scheme}, given by
\begin{equation}
\label{def:BB}
  \Big( GG(v,h) \Big)
  ( u, \tilde{u} )(x)
= \frac{1}{\sqrt{h}}\Bigg(g\left(x, v(x)+\sqrt{h}\, g(x,v(x))\right)-g(x, v(x))\Bigg) \times u(x) \times \tilde{u}(x)
\end{equation}
for all $v\in V_{\beta}$,  then the scheme \eqref{eq:abstract_scheme} reduces to the concrete scheme \eqref{eq:concrete_scheme}.

Apart from \eqref{condition_AA} and \eqref{condition_AAA}, one also needs to verify Assumption \ref{BB_assum} for \eqref{eq:example_spde}.
\begin{prop}\label{lem:BB1}
Suppose that the bilinear operators $G'(v) G(v)$ and $GG(v,h)$ are given by \eqref{def:con_bilinear} and \eqref{def:BB}, respectively. Then  the operator $GG(v,h)$ and the remainder operator $GG_1(v,h)$ given by \eqref{def:BB1} fulfill the conditions in Assumption \ref{BB_assum}, provided that the conditions \eqref{eq:gj} and \eqref{eq:diffusion} hold.
\end{prop}
{\it Proof.} From \eqref{def:BB} we have
\begin{align}\nonumber
&\|GG(v,h)-GG(w,h)\|_{HS^{(2)}(U_0, H)}^2  = \sum_{i,j \in \mathcal{J}} \mu_i\mu_j \Big\|GG(v,h)(\eta_i,\eta_j)-GG(w,h)(\eta_i,\eta_j)\Big\|_H^2
\\ \nonumber=& \frac{1}{h}\sum_{i,j \in \mathcal{J}} \Bigg[\mu_i\mu_j \int_{(0,1)^d} \Bigg| \Bigg(g\Big(x, v(x)+\sqrt{h}\, g(x,v(x))\Big)-g(x, v(x))\Bigg) \times \eta_i(x) \times \eta_j(x)
\\ \nonumber &\quad\quad\quad\quad\quad
- \Bigg(g\Big(x, w(x)+\sqrt{h}\, g(x,w(x))\Big)-g(x, w(x))\Bigg)\times \eta_i(x) \times \eta_j(x) \Bigg|^2 d x\Bigg]
\\\nonumber \leq & \allowdisplaybreaks  \frac{2}{h} \sum_{i,j \in \mathcal{J}} \Bigg[\mu_i\mu_j \int_{(0,1)^d} \left(\Big| g\Big(x, v(x)+\sqrt{h}\cdot g(x,v(x))\Big)-g\Big(x, w(x)+\sqrt{h}\, g(x,w(x))\Big)\Big|^2\right.
\\ & \quad\quad\quad\quad\quad +\left. \Big|g(x, v(x))-g(x, w(x))\Big|^2 \right)dx \times \|\eta_i\|^2_{C((0,1)^d;\mathbb{R})}\times \|\eta_j\|^2_{C((0,1)^d;\mathbb{R})}\Bigg].
\label{eq:bb}
\end{align}
Using \eqref{eq:diffusion} shows that
\begin{equation}\label{eq:bg}
\begin{split}
|g(x, v(x))-g(x, w(x))|^2 \leq c^2 |v(x)-w(x)|^2,
\end{split}
\end{equation}
and that
\begin{equation}\label{eq:bbg}
\begin{split}
& \Big|g\Big(x, v(x)+\sqrt{h}\, g(x,v(x))\Big)-g\Big(x, w(x)+\sqrt{h}\, g(x,w(x))\Big)\Big|^2
\\
\leq& c^2 \Big|v(x)+\sqrt{h}\, g(x,v(x)) - w(x)-\sqrt{h}\, g(x,w(x))\Big|^2
\\
\leq& 2c^2 |v(x) - w(x)|^2 + 2c^2h \,\Big| g(x,v(x))-g(x,w(x))\Big|^2
\\
\leq& 2c^2(1+Tc^2) |v(x) - w(x)|^2,
\end{split}
\end{equation}
where we also used the fact that $h \leq T$. Inserting \eqref{eq:bg} and \eqref{eq:bbg} into \eqref{eq:bb} yields
\begin{align*}
&\|GG(v,h)-GG(w,h)\|_{HS^{(2)}(U_0, H)}^2
\\
\leq& \frac{2}{h} \sum_{i,j \in \mathcal{J}} \Bigg(\mu_i\mu_j \int_{(0,1)^d} c^2(3+2Tc^2)\,|v(x) - w(x)|^2 dx \times \|\eta_i\|^2_{C((0,1)^d;\mathbb{R})}\times \|\eta_j\|^2_{C((0,1)^d;\mathbb{R})}\Bigg)
\\
\leq& \frac{2c^2(3+2Tc^2)}{h}\sum_{i,j \in \mathcal{J}} \mu_i\mu_j\, \Big(\sup_{j \in \mathcal{J}}\|\eta_j\|^4_{C((0,1)^d;\mathbb{R})}\Big)\,  \|v-w\|_H^2
\\
\leq& \frac{2c^6(3+2Tc^2)(TrQ)^2}{h}\,  \|v-w\|_H^2.
\end{align*}
Now the estimate \eqref{eq:BBLip} in Assumption \ref{BB_assum} is validated on choosing $C_0 > 2c^6(3+2Tc^2)(TrQ)^2$.

For the second estimate \eqref{eq:BB1GC}, we use Taylor's formulae in \eqref{def:BB} and derive that for all $v\in V_{\beta}$ the remainder operator $GG_1(v,h) : U_0 \times U_0 \rightarrow H$ given by \eqref{def:BB1} satisfies
\begin{align}\label{eq:BB1Taylor}
&\Big( GG_1(v,h) \Big)
  ( u, \tilde{u} )(x) = \Big( GG(v,h) \Big)
  ( u, \tilde{u} ) -\Big(G'(v)(G(v)u)\Big)(\tilde{u})
  \nonumber\\
=& \sqrt{h} \int_0^1\Big(\frac{\partial^2 }{\partial y^2}g\Big)\Big(x, v(x)+ r \sqrt{h} g(x,v(x))\Big)\,(1-r)\, g^2(x,v(x)) dr \times u(x) \times \tilde{u}(x).
\end{align}
Combining \eqref{eq:diffusion} and \eqref{eq:BB1Taylor} and taking \eqref{eq:HS} into account show
\begin{align*}
&\|GG_1(v,h)\|_{HS^{(2)}(U_0;H)}^2  = \sum_{i,j \in \mathcal{J}} \mu_i\mu_j \|GG_1(v,h)(\eta_i,\eta_j)\|_H^2
\\ =& h\sum_{i,j \in \mathcal{J}} \Bigg(\mu_i\mu_j \int_{(0,1)^d} \Bigg| \int_0^1\Big(\frac{\partial^2 }{\partial y^2}g\Big)\Big(x, v(x)+ r \sqrt{h} g(x,v(x))\Big)\\
& \quad\quad\quad\quad\quad\quad\times
(1-r)\, g^2(x,v(x))dr \times \eta_i(x) \times \eta_j(x)  \Bigg|^2 d x\Bigg)
\\ \leq & c^2 h \sum_{i,j \in \mathcal{J}} \Bigg(\mu_i\mu_j \int_{(0,1)^d} \Big| g(x,v(x))\Big|^4 dx \times \|\eta_i\|^2_{C((0,1)^d;\mathbb{R})}\times \|\eta_j\|^2_{C((0,1)^d;\mathbb{R})}\Bigg).
\end{align*}
Using the elementary inequality $(|a|+|b|)^p \leq 2^{p-1}(|a|^p +|b|^p)$ for $p \geq 1, a,b \in \mathbb{R}$, and \eqref{eq:diffusion} gives
\begin{align*}
| g(x,v(x))|^4 \leq \Big(|g(x,v(x))-g(x,0)|+ |g(x,0)|\Big)^4
\leq \Big(c|v(x)|+ c\Big)^4
\leq 8c^4 \Big(|v(x)|^4+ 1 \Big).
\end{align*}
Thus for all $v \in V_{\beta}$,
\begin{align*}
\|GG_1(v,h)\|_{HS^{(2)}(U_0;H)}^2 \leq& 8c^6 h \sum_{i,j \in \mathcal{J}} \mu_i\mu_j\, \Big(\sup_{j \in \mathcal{J}}\|\eta_j\|^4_{C((0,1)^d;\mathbb{R})}\Big)\, \Big(\|v\|_{L^4((0,1)^d;\mathbb{R})}^4+1 \Big)
\\ \leq & 8c^{10}(TrQ)^2 h \, (\|v\|_{V_{\beta}}^4+1),
\end{align*}
where \eqref{eq:gj} and the fact were used that  $V_{\beta} \subset L^4((0,1)^d;\mathbb{R})$ continuously for $\beta = \frac{d}{5}$ by Sobolev embedding theorem. The proof of Proposition \ref{lem:BB1} is complete.
$\square$

So far, all conditions in Assumption \ref{semigroup}-\ref{BB_assum} have been verified except conditions \eqref{condition_AA} and \eqref{condition_AAA} in Assumption \ref{diffusion}. We mention that condition \eqref{condition_AA} originally comes from \cite{JR10a}, where \eqref{condition_AA} is a key ingredient to promise higher spatial and temporal regularity of mild solution. The condition \eqref{condition_AAA} is needed to estimate approximation error due to truncation of the expansion of the noise $W_t$. In the last section, we will obtain these two conditions for two concrete examples.

\section{Proof of Theorem \ref{thm:mainresult}}\label{sec:proof}

First of all, we rewrite the numerical solution~\eqref{eq:abstract_scheme} in the following form
\begin{align}\label{eq:defYNMK}
  Y_{ m }^{ N,M,K }
&=
  S( m h ) P_N\!\left( \xi \right)
  +
  P_N\!\left(
    \sum_{ l=0 }^{ m-1 }
    \int_{ l h }^{ (l+1)h }
    S\Big((m-l)h\Big)
    F\!\left( Y_{ l }^{ N,M,K } \right)
    ds
  \right)
\nonumber
\\&\quad+
  P_N\!\left(
    \sum_{ l=0 }^{ m-1 }
    \int_{ l h }^{ (l+1)h }
    S\Big((m-l)h\Big)
    G\!\left( Y_{ l }^{ N,M,K } \right)
    dW_s^K
  \right)
\\&\quad+
\nonumber
  \frac{1}{2}P_N\!\left(
    \sum_{ l=0 }^{ m-1 }S\Big((m-l)h\Big) GG(Y_l^{N,M,K},h)(\Delta W^{M,K}_l, \Delta W^{M,K}_l)
  \right)
\\ & \quad -
\nonumber
  \frac{h}{2}P_N\!\Bigg(
    \sum_{ l=0 }^{ m-1 }S\Big((m-l)h\Big) \sum_{\begin{subarray}{ll} j\in \mathcal{J}_K\\ \mu_j \neq 0 \end{subarray}} \mu_j GG(Y_l^{N,M,K},h)(\eta_j, \eta_j)
  \Bigg).
\end{align}
Likewise, the numerical solution \eqref{eq:Mscheme} satisfies
\begin{align}\label{eq:MilYP}
  \tilde{Y}_{ m }^{ N,M,K }
&=
  S( m h ) P_N\!\left( \xi \right)
  +
  P_N\!\left(
    \sum_{ l=0 }^{ m-1 }
    \int_{ l h }^{ (l+1)h }
    S\Big((m-l)h\Big)
    F\!\left( \tilde{Y}_{ l }^{ N,M,K } \right)
    ds
  \right)
\nonumber
\\&\quad+
  P_N\!\left(
    \sum_{ l=0 }^{ m-1 }
    \int_{ l h }^{ (l+1)h }
    S\Big((m-l)h\Big)
    G\!\left( \tilde{Y}_{ l }^{ N,M,K } \right)
    dW_s^K
  \right)
\\&\quad+
\nonumber
  \frac{1}{2}P_N\!\left(
    \sum_{ l=0 }^{ m-1 }S\Big((m-l)h\Big) G'(\tilde{Y}_{ l }^{ N,M,K })\Big(G(\tilde{Y}_{ l }^{ N,M,K })\Delta W^{M,K}_l\Big) \Big(\Delta W^{M,K}_l\Big)
  \right)
\\ & \quad -
\nonumber
  \frac{h}{2}P_N\!\Bigg(
    \sum_{ l=0 }^{ m-1 }S\Big((m-l)h\Big) \sum_{\begin{subarray}{ll} j\in \mathcal{J}_K\\ \mu_j \neq 0 \end{subarray}} \mu_j G'(\tilde{Y}_{ l }^{ N,M,K })\Big(G(\tilde{Y}_{ l }^{ N,M,K })\eta_j\Big) \eta_j
  \Bigg).
\end{align}
Moreover, the exact mild solution of
the SPDE~\eqref{eq:abstract_spde} can be rewritten as
\begin{align}\label{eq:exactrep}
  X_{ m h }
&=
  S( m h ) \xi
  +
  \int_0^{ m h }
  S(m h - s)
  F\!\left( X_s \right)
  ds
  +
  \int_0^{ m h }
  S( m h - s)
  G\!\left( X_s \right)
  dW_s
\\&=
\nonumber
  S( m h ) \xi
  +
  \sum_{ l=0 }^{ m-1 }
  \int_{ l h }^{ (l+1)h }
  S(m h - s)
  F\!\left( X_s \right)
  ds
  +
  \sum_{ l=0 }^{ m-1 }
  \int_{ l h }^{ (l+1)h }
  S(m h - s)
  G\!\left( X_s \right)
  dW_s,
\end{align}
and thus
\begin{align}\label{eq:PNX}
  P_N\!\left( X_{ m h } \right)
&=
  S( m h ) P_N( \xi )
  +
  P_N\!\left(
  \sum_{ l=0 }^{ m-1 }
  \int_{ l h }^{ (l+1)h }
  S( m h - s)
  F\!\left( X_s \right)
  ds
  \right)
\nonumber
\\&\quad+
  P_N\!\left(
  \sum_{ l=0 }^{ m-1 }
  \int_{ l h }^{ (l+1)h }
  S(m h-s)
  G\!\left( X_s \right)
  dW_s \right).
\end{align}

To estimate
$ \mathbb{E}
\left\| X_{ m h } - Y_m^{ N,M,K }
\right\|_H^2 $, we need two auxiliary processes $Z_m^{N, M, K}$ and $ \tilde{Z}_m^{N, M, K} $. Define $Z_m^{N, M, K}$ by
\begin{align}\label{eq:defZ}
  Z_{ m }^{ N, M, K }
&:=
  S( m h ) P_N\!\left( \xi \right)
  +
  P_N\!\left(
    \sum_{ l=0 }^{ m-1 }
    \int_{ l h }^{ (l+1)h }
    S\Big((m-l)h\Big)
    F\!\left( X_{ lh } \right)
    ds
  \right)
\nonumber
\\&\quad+
  P_N\!\left(
    \sum_{ l=0 }^{ m-1 }
    \int_{ l h }^{ (l+1)h }
    S\Big((m-l)h \Big)
    G\!\left( X_{ lh } \right)
    dW_s^K
  \right)
\\&\quad+
\nonumber
  \frac{1}{2}P_N\!\left(
    \sum_{ l=0 }^{ m-1 }S\Big((m-l)h\Big) GG(X_{lh},h)(\Delta W^{M,K}_l, \Delta W^{M,K}_l)
  \right)
\\ & \quad -
\nonumber
  \frac{h}{2}P_N\!\Bigg(
    \sum_{ l=0 }^{ m-1 }S\Big((m-l)h\Big) \sum_{\begin{subarray}{ll} j\in \mathcal{J}_K\\ \mu_j \neq 0 \end{subarray}} \mu_j GG(X_{lh},h)(\eta_j, \eta_j)
  \Bigg).
\end{align}
Note that \eqref{eq:defYNMK} coincides with \eqref{eq:defZ} with $Y_l^{N,M,K}$ replaced by $X_{lh}$. Similarly, replacing $\tilde{Y}_{ l }^{ N,M,K }$ in \eqref{eq:MilYP} by $X_{lh}$ we introduce the process $ \tilde{Z}_m^{N, M, K} $, given by
\begin{align}\label{eq:defZt}
  \tilde{Z}_{ m }^{ N, M, K }
&:=
  S( m h ) P_N\!\left( \xi \right)
  +
  P_N\!\left(
    \sum_{ l=0 }^{ m-1 }
    \int_{ l h }^{ (l+1)h }
    S\Big((m-l)h\Big)
    F\!\left( X_{ lh } \right)
    ds
  \right)
\nonumber
\\&\quad+
  P_N\!\left(
    \sum_{ l=0 }^{ m-1 }
    \int_{ l h }^{ (l+1)h }
    S\Big((m-l)h\Big)
    G\!\left( X_{ lh } \right)
    dW_s^K
  \right)
\\&\quad+
\nonumber
  \frac{1}{2}P_N\!\left(
    \sum_{ l=0 }^{ m-1 }S\Big((m-l)h\Big) G'(X_{lh})\Big(G(X_{lh})\Delta W^{M,K}_l\Big) \Big(\Delta W^{M,K}_l\Big)
  \right)
\\ & \quad -
\nonumber
  \frac{h}{2}P_N\!\Bigg(
    \sum_{ l=0 }^{ m-1 }S\Big((m-l)h\Big) \sum_{\begin{subarray}{ll} j\in \mathcal{J}_K\\ \mu_j \neq 0 \end{subarray}} \mu_j G'(X_{lh})\Big(G(X_{lh})\eta_j\Big) \eta_j
  \Bigg).
\end{align}

Armed with these notations, now we start the proof. Employing the elementary inequality
$( a_1 + a_2 + a_3)^2
\leq
  3 ( |a_1|^2 + |a_2|^2 + |a_3|^2)
$, $a_1, a_2, a_3 \in \mathbb{R}$ shows for $m= 0,1,...,M$
\begin{align}\label{eq:beginproof}
&
  \mathbb{E}\left\|
    X_{ m h } - Y_m^{ N,M,K }
  \right\|_H^2
\\\leq &
\nonumber
  3
  \mathbb{E}\left\|
    X_{ m h } - P_N\!\left( X_{ m h } \right)
  \right\|_H^2
  +
  3 \mathbb{E}\left\|
    P_N\!\left( X_{ m h } \right)
    - Z_m^{ N, M,K }
  \right\|_H^2
  +
  3 \mathbb{E}\left\|
    Z_m^{ N, M,K } - Y_m^{ N,M,K }
  \right\|_H^2.
\end{align}
Below we will estimate the three terms in \eqref{eq:beginproof} step by step. First let $ R >0 $ be a real constant
such that
\begin{equation}\label{eq:R}
\begin{split}
  &\left\| F'\!\left( v \right) \right\|_{ L( H ) }
  \leq R,
  \qquad
  \left\| G'\!\left( v \right)
  \right\|_{ L( H, HS( U_0, H ) ) }
  \leq R ,\\&
  \mathbb{E}\left\| \left( -A \right)^{ \gamma }
    X_t \right\|_{ H }^2
  =
  \mathbb{E}\left\| X_t
  \right\|_{ V_{\gamma} }^2
  \leq R,
  \;\;
  \mathbb{E}\left\| X_t
  \right\|^4_{V_{\beta}}
  \leq
  R
\end{split}
\end{equation}
for all $ v \in V_{\beta} $
and all $ t\in [0,T] $.
Due to Assumptions~\ref{semigroup}-\ref{initial} in Section~\ref{sec:settings}
and Proposition~\ref{lem:existence},
such a real constant exists.

For the spatial discretization error
$ \mathbb{E}\left\|
    X_{ m h }
    - P_N\!\left(
      X_{ m h }
    \right)
  \right\|_H^2 $,
using \eqref{eq:R} we derive
\begin{align}\nonumber
  \mathbb{E}\left\|
    X_{ m h }
    - P_N\!\left(
      X_{ m h }
    \right)
  \right\|_H^2
&=
  \mathbb{E}\left\|
    \left( I - P_N \right) X_{ m h }
  \right\|_H^2
\nonumber
\\&\leq
  \left\|
    \left( - A \right)^{ - \gamma }
    \left( I - P_N \right)
  \right\|_{ L(H) }^2 \times
  \mathbb{E}\left\|X_{ m h }
  \right\|_{V_{\gamma}}^2
\nonumber
\\& \leq R   \left\|
    \left( - A \right)^{ - \gamma }
    \left( I - P_N \right)
  \right\|_{ L(H) }^2 = R \left(
    \inf_{
      i \in \mathcal{I} \backslash
      \mathcal{I}_N
    }
    \lambda_i
  \right)^{ -2\gamma }.
\label{eq:spatialerror}
\end{align}


To estimate $\mathbb{E} \|P_N(X_{mh}) - Z_{ m }^{ N, M, K }\|^2_H$, we need the estimate $\mathbb{E} \|P_N(X_{mh}) - \tilde{Z}_{ m }^{ N, M, K }\|^2_H$.
\begin{lemma}\label{lem:Mil_Loc_Error}
Under Assumptions \ref{semigroup}-\ref{initial} and \eqref{eq:trQ}, there exists a constant $C_1$, independent of $h$, such that
\begin{align}\label{eq:Milstein_error}
\sup_{0\leq m \leq M}\mathbb{E} \Big\|P_N(X_{mh}) - \tilde{Z}_{ m }^{ N, M, K }\Big\|^2_H
      \leq      C_1 \left( \bigg(
      \sup_{ j \in \mathcal{J} \backslash \mathcal{J}_K }
      \mu_j
      \bigg)^{ 2\alpha }
    +
      M^{ - \min\left( 4 \left( \gamma - \beta \right) ,
      2\gamma \right) } \right).
\end{align}
\end{lemma}
{\it Proof.} To establish the convergence result for the Milstein type scheme \eqref{eq:Mscheme},
 \eqref{eq:Milstein_error}
has been obtained in \cite{JR10b} (see Section 5 in \cite{JR10b} for the details). $\square$

\begin{lemma} Suppose that Assumptions~\ref{semigroup}-\ref{BB_assum} and the condition \eqref{eq:trQ} are fulfilled.
Then there exists a constant  $C_2$, independent of $h$, such that
\begin{equation}\label{eq:LemTemp_error}
\sup_{0\leq m \leq M} \mathbb{E} \Big\|P_N(X_{mh}) - Z_{ m }^{ N, M, K }\Big\|^2_H \leq C_2\bigg(
      \sup_{ j \in \mathcal{J} \backslash \mathcal{J}_K }
      \mu_j
      \bigg)^{ 2\alpha }
    +
      \frac{C_2}{M^{\min\left( 4 \left( \gamma - \beta \right) ,
      2\gamma \right)}}.
\end{equation}
\end{lemma}

{\it Proof.} Using \eqref{def:BB1}, we derive from \eqref{eq:defZ} and \eqref{eq:defZt} that
\begin{align}\label{eq:ZZT}
Z_{ m }^{ N, M, K } =& \tilde{Z}_{ m }^{ N, M, K } + \frac{1}{2}P_N\!\left(
    \sum_{ l=0 }^{ m-1 }S\Big((m-l)h\Big) GG_1(X_{lh},h)(\Delta W^{M,K}_l, \Delta W^{M,K}_l)
  \right)
\\ & -
\nonumber
  \frac{h}{2}P_N\!\Bigg(
    \sum_{ l=0 }^{ m-1 }S\Big((m-l)h\Big) \sum_{\begin{subarray}{ll} j\in \mathcal{J}_K\\ \mu_j \neq 0 \end{subarray}} \mu_j GG_1(X_{lh},h)(\eta_j, \eta_j)
  \Bigg).
\end{align}
Therefore
\begin{align}\label{eq:PXZ}
P_N(X_{mh}) - Z_{ m }^{ N, M, K } = & P_N(X_{mh}) - \tilde{Z}_{ m }^{ N, M, K }
\nonumber \\ &- \frac{1}{2}P_N\!\left(
    \sum_{ l=0 }^{ m-1 }S\Big((m-l)h\Big) GG_1(X_{lh},h)(\Delta W^{M,K}_l, \Delta W^{M,K}_l)
  \right)
\nonumber\\ & +
\nonumber
  \frac{h}{2}P_N\!\Bigg(
    \sum_{ l=0 }^{ m-1 }S\Big((m-l)h\Big) \sum_{\begin{subarray}{ll} j\in \mathcal{J}_K\\ \mu_j \neq 0 \end{subarray}} \mu_j GG_1(X_{lh},h)(\eta_j, \eta_j)
  \Bigg),
\end{align}
and thus using the elementary inequality $(a_1+a_2)^2 \leq 2(|a_1|^2+|a_2|^2), a_1, a_2 \in \mathbb{R}$  yields
\begin{equation}\label{eq:Temp_error}
\begin{split}
\mathbb{E} \Big\|P_N(X_{mh}) - Z_{ m }^{ N, M, K }\Big\|^2_H \leq&  2\mathbb{E} \Big\|P_N(X_{mh}) - \tilde{Z}_{ m }^{ N, M, K }\Big\|^2_H
\\ & +  \frac{1}{2} \mathbb{E} \Bigg\|
    \sum_{ l=0 }^{ m-1 }S\Big((m-l)h\Big) GG_1(X_{lh},h)\Big(\Delta W^{M,K}_l, \Delta W^{M,K}_l\Big)
\\
 &-  h\sum_{ l=0 }^{ m-1 }S\Big((m-l)h\Big) \sum_{\begin{subarray}{ll} j\in \mathcal{J}_K\\ \mu_j \neq 0 \end{subarray}} \mu_j GG_1(X_{lh},h)(\eta_j, \eta_j) \Bigg\|^2_H
\\
:=& 2\mathbb{E} \Big\|P_N(X_{mh}) - \tilde{Z}_{ m }^{ N, M, K }\Big\|^2_H + J_1,
\end{split}
\end{equation}
where the fact that $\|P_N v\|_H \leq \|v\|_H$ was also used. Due to \eqref{eq:Milstein_error} in Lemma \ref{lem:Mil_Loc_Error}, it remains to estimate $J_1$.
Inserting the representation \eqref{Wiener_Incr} and using bilinearity of the operator $GG_1$ give
\begin{align}\nonumber
  J_1
=&
  \frac{1}{2}\mathbb{E}\Big\|
    \sum_{ l=0 }^{ m-1 }S\Big((m-l)h\Big) \Big(GG_1(X_{lh},h)\Big)\Big( \sum_{\begin{subarray}{ll} i\in \mathcal{J}_K\\ \mu_i \neq 0 \end{subarray}} \sqrt{\mu_i} \Delta \beta_l^i \eta_i, \sum_{\begin{subarray}{ll} j\in \mathcal{J}_K\\ \mu_j \neq 0 \end{subarray}} \sqrt{\mu_j} \Delta \beta_l^j \eta_j\Big)
\\ \nonumber & \quad-h
    \sum_{ l=0 }^{ m-1 }S\Big((m-l)h\Big) \sum_{\begin{subarray}{ll} i,j\in \mathcal{J}_K\\ \mu_i,\mu_j \neq 0 \end{subarray}} \sqrt{\mu_i} \sqrt{\mu_j}\delta_{ij} \Big(GG_1(X_{lh},h)\Big)(\eta_i, \eta_j)
  \Bigg\|_{ H }^2
\\=&
  \frac{1}{2}\mathbb{E}\Bigg\|
    \sum_{ l=0 }^{ m-1 }\sum_{\begin{subarray}{ll} i,j\in \mathcal{J}_K\\ \mu_i,\mu_j \neq 0 \end{subarray}}\sqrt{\mu_i} \sqrt{\mu_j} S\Big((m-l)h\Big) \Big(GG_1(X_{lh},h)\Big)( \eta_i, \eta_j)\times (\Delta \beta_l^i\Delta \beta_l^j- \delta_{ij}h)  \Bigg\|_{ H }^2,
\label{eq:J_3_2}
\end{align}
where $\delta_{ij} = 1$ for $i=j$ and $\delta_{ij} = 0$ for $i\neq j$. For simplicity of notation, we denote
$$
\chi_l^{i,j} = \sqrt{\mu_i} \sqrt{\mu_j} S\Big((m-l)h\Big) \Big(GG_1(X_{lh},h)\Big)( \eta_i, \eta_j).
$$
Then we can rewrite \eqref{eq:J_3_2} as
\begin{align*}
  J_1
=&
  \frac{1}{2}\mathbb{E}\Bigg\langle
    \sum_{ l_1=0 }^{ m-1 }\sum_{\begin{subarray}{ll} i_1,j_1\in \mathcal{J}_K\\ \mu_{i_1},\mu_{j_1} \neq 0 \end{subarray}} \chi_{l_1}^{i_1,j_1}\, (\Delta \beta_{l_1}^{i_1}\Delta \beta_{l_1}^{j_1}- \delta_{i_1j_1}h), \sum_{ l_2=0 }^{ m-1 }\sum_{\begin{subarray}{ll} i_2,j_2\in \mathcal{J}_K\\ \mu_{i_2},\mu_{j_2} \neq 0 \end{subarray}} \chi_{l_2}^{i_2,j_2} \, (\Delta \beta_{l_2}^{i_2}\Delta \beta_{l_2}^{j_2}- \delta_{i_2j_2}h)\Bigg\rangle_H
\\ =& \frac{1}{2} \sum_{ l_1=0 }^{ m-1 } \sum_{ l_2=0 }^{ m-1 } \sum_{\begin{subarray}{ll} i_1,j_1\in \mathcal{J}_K\\ \mu_{i_1},\mu_{j_1} \neq 0 \end{subarray}} \sum_{\begin{subarray}{ll} i_2,j_2\in \mathcal{J}_K\\ \mu_{i_2},\mu_{j_2} \neq 0 \end{subarray}} \mathbb{E}\Bigg(\Big\langle
     \chi_{l_1}^{i_1,j_1},  \chi_{l_2}^{i_2,j_2} \Big\rangle_H \, (\Delta \beta_{l_1}^{i_1}\Delta \beta_{l_1}^{j_1}- \delta_{i_1j_1}h)\, (\Delta \beta_{l_2}^{i_2}\Delta \beta_{l_2}^{j_2}- \delta_{i_2j_2}h)\Bigg).
\end{align*}
In the case that $l_1 \neq l_2$, without loss of generality we set $l_1 < l_2$. Using the fact that $(\Delta \beta_{l_2}^i)_{i \in \mathcal{J}_K, \mu_i \neq 0}$ are independent of $\mathcal{F}_{t_{l_2}} \supset \mathcal{F}_{t_{l_1}}$, $\mathbb{E} (\Delta \beta_{l}^i) =0$ and $\chi_l^{i,j} \in \mathcal{F}_{t_l}$ shows that
\begin{align*}
&\mathbb{E}\left[\Big\langle
     \chi_{l_1}^{i_1,j_1},  \chi_{l_2}^{i_2,j_2} \Big\rangle_H \, (\Delta \beta_{l_1}^{i_1}\Delta \beta_{l_1}^{j_1}- \delta_{i_1j_1}h)\, (\Delta \beta_{l_2}^{i_2}\Delta \beta_{l_2}^{j_2}- \delta_{i_2j_2}h)\right]
\\ =&  \mathbb{E} \left[\Big\langle
     \chi_{l_1}^{i_1,j_1},  \chi_{l_2}^{i_2,j_2} \Big\rangle_H \, (\Delta \beta_{l_1}^{i_1}\Delta \beta_{l_1}^{j_1}- \delta_{i_1j_1}h)\right]\, \mathbb{E}\Big( \Delta \beta_{l_2}^{i_2}\Delta \beta_{l_2}^{j_2}- \delta_{i_2j_2}h\Big)
\\ = & 0.
\end{align*}
Here the last step follows by the obvious fact that $\mathbb{E}\left( \Delta \beta_{l_2}^{i_2}\Delta \beta_{l_2}^{j_2}- \delta_{i_2j_2}h\right) = 0$.
In the case that $l_1 = l_2 = l$ but that $\{i_1,j_1\} \neq \{i_2, j_2\}$. Using the mutual independence of $\Delta \beta_l^i, \Delta \beta_l^j, i \neq j$ and the fact that $(\Delta \beta_{l}^i)_{i \in \mathcal{J}_K, \mu_i \neq 0}$ are independent of $\mathcal{F}_{t_{l}}$, $\mathbb{E} (\Delta \beta_{l}^i) =0$  and $\chi_l^{i,j} \in \mathcal{F}_{t_l}$ shows that
\begin{align*}
&\mathbb{E}\left[\Big\langle
     \chi_{l_1}^{i_1,j_1},  \chi_{l_2}^{i_2,j_2} \Big\rangle_H \, (\Delta \beta_{l_1}^{i_1}\Delta \beta_{l_1}^{j_1}- \delta_{i_1j_1}h)\, (\Delta \beta_{l_2}^{i_2}\Delta \beta_{l_2}^{j_2}- \delta_{i_2j_2}h)\right]
\\ =&  \mathbb{E}\Big\langle
     \chi_{l}^{i_1,j_1},  \chi_{l}^{i_2,j_2} \Big\rangle_H\, \mathbb{E}\Big( (\Delta \beta_{l}^{i_1}\Delta \beta_{l}^{j_1}- \delta_{i_1j_1}h)\,(\Delta \beta_{l}^{i_2}\Delta \beta_{l}^{j_2}- \delta_{i_2j_2}h)\Big)
\\ = & 0.
\end{align*}
Here the last step follows since
$$\mathbb{E}\left( (\Delta \beta_{l}^{i_1}\Delta \beta_{l}^{j_1}- \delta_{i_1j_1}h)\,(\Delta \beta_{l}^{i_2}\Delta \beta_{l}^{j_2}- \delta_{i_2j_2}h)\right) = 0
$$
in the case that $\{i_1,j_1\} \neq \{i_2, j_2\}$.
Hence, using the fact that $\chi_l^{i,j} \in \mathcal{F}_{t_l}$ and that $(\Delta \beta_{l}^i)_{i \in \mathcal{J}_K, \mu_i \neq 0}$ are independent of $\mathcal{F}_{t_{l}}$ we have
\begin{align}\label{eq:J_3}
  J_1
=&
  \frac{1}{2}\sum_{ l=0 }^{ m-1 }\sum_{\begin{subarray}{ll} i,j\in \mathcal{J}_K\\ \mu_i,\mu_j \neq 0 \end{subarray}} \mathbb{E}\Big\|
    \chi_l^{i,j} \, (\Delta \beta_l^i\Delta \beta_l^j- \delta_{ij}h)  \Big\|_{ H }^2.
\nonumber\\ \leq&
  \sum_{ l=0 }^{ m-1 } \sum_{\begin{subarray}{ll} i,j\in \mathcal{J}_K\\ \mu_i,\mu_j \neq 0 \end{subarray}}\mathbb{E}\Big\|
      \chi_l^{i,j}\, (\Delta \beta_l^i\Delta \beta_l^j)  \Big\|_{ H }^2
 + h^2  \sum_{ l=0 }^{ m-1 } \mathbb{E} \Bigg(\sum_{\begin{subarray}{ll} i,j\in \mathcal{J}_K\\ \mu_i,\mu_j \neq 0 \end{subarray}}\Big\|
      \chi_l^{i,j}  \Big\|_{ H }^2\Bigg)
\nonumber\\ =&
  \sum_{ l=0 }^{ m-1 } \sum_{\begin{subarray}{ll} i,j\in \mathcal{J}_K\\ \mu_i,\mu_j \neq 0 \end{subarray}}\Bigg(\mathbb{E}\Big\| S\Big((m-l)h\Big)
       \Big(GG_1(X_{lh},h)\Big)( \sqrt{\mu_i}\eta_i, \sqrt{\mu_j} \eta_j) \Big\|_{ H }^2 \, \mathbb{E}|\Delta \beta_l^i\Delta \beta_l^j|^2\Bigg)
\nonumber\\ & + h^2  \sum_{ l=0 }^{ m-1 } \mathbb{E} \Bigg(\sum_{\begin{subarray}{ll} i,j\in \mathcal{J}_K\\ \mu_i,\mu_j \neq 0 \end{subarray}}\Big\|S\Big((m-l)h\Big)
      \Big(GG_1(X_{lh},h)\Big)( \sqrt{\mu_i}\eta_i, \sqrt{\mu_j} \eta_j)  \Big\|_{ H }^2\Bigg).\nonumber
\end{align}
Due to the fact that $\|S(t) v\|_H \leq \|v\|_H, t \geq 0$ for all $v \in H$ and that $\mathbb{E} |\Delta \beta_l^i\Delta \beta_l^j|^2 \leq 3h^2$ for all $i,j \in \mathcal{J}_K$, and using \eqref{eq:HS}, \eqref{eq:BB1GC} in Assumption \ref{BB_assum} and \eqref{eq:R} we derive that
\begin{align}\nonumber
J_1\leq &
3 h^2\sum_{ l=0 }^{ m-1 } \mathbb{E} \Bigg(\sum_{\begin{subarray}{ll} i,j\in \mathcal{J}_K \\ \mu_i,\mu_j \neq 0 \end{subarray}}\Bigg\|
      \Big(GG_1(X_{lh},h)\Big)( \sqrt{\mu_i}\eta_i, \sqrt{\mu_j} \eta_j)  \Bigg\|_{ H }^2\Bigg)
\nonumber \\ &
      + h^2 \sum_{l=0}^{m-1} \mathbb{E}\Big\|GG_1(X_{lh},h)\Big\|_{HS^{(2)}(U_0, H)}^2
\nonumber \\
      \leq& 4h^2 \sum_{l=0}^{m-1} \mathbb{E}\Big\|GG_1(X_{lh},h)\Big\|_{HS^{(2)}(U_0, H)}^2
\nonumber \\
       \leq&  4C_0h^3 \sum_{l=0}^{m-1} (1+ \mathbb{E}\| X_{lh}\|_{V_{\beta}}^4)
\nonumber \\
       \leq & \frac{4C_0T^3(1+R)}{M^2}.
\end{align}
Plugging \eqref{eq:Milstein_error} and the preceding estimate into \eqref{eq:Temp_error} gives
\begin{align}
\mathbb{E} \Big\|P_N(X_{mh}) - Z_{ m }^{ N, M, K }\Big\|^2_H \leq 2C_1\bigg(
      \sup_{ j \in \mathcal{J} \backslash \mathcal{J}_K }
      \mu_j
      \bigg)^{ 2\alpha }
    +
      \frac{2C_1+4C_0T^3(1+R)}{M^{\min\left( 4 \left( \gamma - \beta \right) ,
      2\gamma \right)}}.
\end{align}
Consequently, \eqref{eq:LemTemp_error} is derived on choosing $C_2 = 2C_1+4cT^3(1+R)$.
$\square$

\begin{lemma} Suppose that all conditions in Assumptions~\ref{semigroup}-\ref{BB_assum} are fulfilled. Then there exists a constant $C_3$, independent of $h$, such that
\begin{equation}\label{eq:LemLips}
\mathbb{E}\left\|
    Z_m^{ N, M, K }
    -
    Y_m^{ N, M, K }
  \right\|_{ H }^2 \leq C_3h \sum_{l=0}^{m-1} \mathbb{E}\Big\| X_{lh}-Y_l^{ N,M,K } \Big\|_H^2
\end{equation}
holds for all $m =0,1,\dots, M$.
\end{lemma}

{\it Proof.} Using the elementary inequality $( a_1 + a_2 + a_3)^2
\leq
  3 ( |a_1|^2 + |a_2|^2 + |a_3|^2), a_1, a_2, a_3 \in\mathbb{R}
$ gives
\begin{align}\nonumber
  &\mathbb{E}\left\|
    Z_m^{ N, M, K }
    -
    Y_m^{ N, M, K }
  \right\|_{ H }^2
\\ \leq&
  3  \mathbb{E}\left\|
    \sum_{ l=0 }^{ m - 1 }
    \int_{ l h }^{ (l+1)h }
    S\Big((m-l)h\Big)
    \left(
      F\!\left( X_{ l h } \right)
      -
      F\!( Y_l^{ N,M,K })
    \right) ds
  \right\|_{ H }^2
\nonumber\\&+
 3  \mathbb{E}\left\|
    \sum_{ l=0 }^{ m - 1 }
    \int_{ l h }^{ (l+1)h }
    S\Big((m-l)h\Big)
    \left(
      G\!\left( X_{ l h } \right)
      -
      G\!( Y_l^{ N,M,K })
    \right) dW_s^K
  \right\|_{ H }^2
\nonumber\\&+
  \frac{3}{4}\mathbb{E}\Bigg\|
    \sum_{ l=0 }^{ m-1 }S\Big((m-l)h\Big) \Big(GG(X_{lh},h)-GG(Y_l^{ N,M,K },h)\Big)(\Delta W^{M,K}_l, \Delta W^{M,K}_l)
\nonumber\\ &-h
    \sum_{ l=0 }^{ m-1 }S\Big((m-l)h\Big) \sum_{\begin{subarray}{ll} j\in \mathcal{J}_K\\ \mu_j \neq 0 \end{subarray}} \mu_j \Big(GG(X_{lh},h)-GG(Y_l^{ N,M,K },h)\Big)(\eta_j, \eta_j)
  \Bigg\|_{ H }^2
\nonumber\\:=& J_{2}+J_{3}+J_{4}, \label{eq:Lips}
\end{align}
where we also used the fact that $\|P_N v\|_H \leq \|v\|_H$.
For $J_2$, one can derive that
\begin{align}\label{eq:J2}
 J_2 \leq& 3 M h^2 \left(
    \sum_{ l=0 }^{ m - 1 }
    \mathbb{E}\left\|
      S\Big((m-l)h\Big)
      \left(
        F\!\left( X_{ l h } \right)
        -
        F\!\left( Y_l^{ N,M,K } \right)
      \right)
    \right\|_{ H }^2
  \right)
\nonumber \\
\leq & 3Th \Bigg(\sum_{ l=0 }^{ m - 1 }\mathbb{E}\left\|
        F\!\left( X_{ l h } \right)
        -
        F\!\left( Y_l^{ N,M,K } \right)
    \right\|_{ H }^2\Bigg)
\nonumber \\
\leq & 3TR^2h \sum_{l=0}^{m-1} \mathbb{E}\Big\|X_{lh}-Y_l^{ N,M,K }\Big\|_H^2,
\end{align}
where \eqref{eq:R} and the fact that $\|S(t)\|_{L(H)} \leq 1, t \geq 0$ were used.
For $J_3$, one can derive that
\begin{align}\label{eq:J3}
 J_3 \leq & 3 \sum_{ l=0 }^{ m - 1 }
    \int_{ l h }^{ (l+1)h }\mathbb{E}
    \left\| S\Big((m-l)h\Big)
    \left(
      G\!\left( X_{ l h } \right)
      -
      G\!\left( Y_l^{ N,M,K } \right)
    \right) \right\|_{ HS(U_0,H) }^2 ds
\nonumber \\
\leq & 3R^2h \sum_{l=0}^{m-1} \mathbb{E}\Big\|X_{lh}-Y_l^{ N,M,K }\Big\|_H^2,
\end{align}
where the isometry property, \eqref{eq:R} and the fact that $\|S(t)\|_{L(H)} \leq 1, t \geq 0$ were again used.

Now it remains to estimate $J_{4}$. In a similar way as estimating $J_1$ and using the condition \eqref{eq:BBLip} in Assumption \ref{BB_assum} and \eqref{eq:R},  one can obtain that
\begin{align} \label{eq:J4}
J_4 \leq&  6h^2 \sum_{l=0}^{m-1} \mathbb{E}\|GG(X_{lh},h)-GG(Y_l^{N,M,K},h)\|_{HS^{(2)}(U_0, H)}^2
\nonumber \\
      \leq&  6C_0h \sum_{l=0}^{m-1} \mathbb{E}\|X_{lh} - Y_l^{N,M,K}\|_{H}^2.
\end{align}
Now, inserting \eqref{eq:J2}, \eqref{eq:J3} and \eqref{eq:J4} into \eqref{eq:Lips} gives the desired estimate \eqref{eq:LemLips}. $\square$

Now we return to \eqref{eq:beginproof}.  With the estimates \eqref{eq:spatialerror}, \eqref{eq:LemTemp_error} and \eqref{eq:LemLips} at hand, we derive from \eqref{eq:beginproof} that
\begin{align}\label{eq:endproof}
  \mathbb{E}\left\|
    X_{ m h } - Y_m^{ N,M,K }
  \right\|_H^2 \leq & 3C_3 \, h \,\sum_{l=0}^{m-1} \mathbb{E}\Big\| X_{lh}-Y_l^{ N,M,K } \Big\|_H^2
\\& +
\nonumber
  3R \left(
    \inf_{
      i \in \mathcal{I} \backslash
      \mathcal{I}_N
    }
    \lambda_i
  \right)^{ -2\gamma }
  +
  3C_2\bigg(
      \sup_{ j \in \mathcal{J} \backslash \mathcal{J}_K }
      \mu_j
      \bigg)^{ 2\alpha }
    +
      \frac{3C_2}{M^{\min\left( 4 \left( \gamma - \beta \right) ,
      2\gamma \right)}}.
\end{align}
Finally, Gronwall's lemma gives the main result \eqref{eq:mainresult} and the proof of Theorem \ref{thm:mainresult} is complete.

\section{Numerical experiments}\label{sec:Numer-Experi}

In this section we will first illustrate how to implement the scheme introduced in this work and then present some numerical results to support our theoretical assertions.

\subsection{Implementation}

For simplicity of notations, here we only consider the new scheme \eqref{eq:concrete_scheme} for one dimensional space case \eqref{eq:spde} and one can adapt the following implementation to handle multi-dimensional space case and other schemes. Using the notation $\zeta_m : (0,1)\rightarrow \mathbb{R}$ given by
\begin{equation}\label{eq:zeta}
\begin{split}
\zeta_m (x) = & Y_m^{N,M,K}(x) + h \, f(x, Y_m^{N,M,K}(x)) + g\left(x, Y_m^{N,M,K}(x)\right)\times \Delta W_m^{M,K}(x)
\\  & +\frac{1}{2\sqrt{h}} \left[ g\Big(x, Y_m^{N,M,K}(x) + \sqrt{h}\,g(x, Y_m^{N,M,K}(x))\Big)-g(x, Y_m^{N,M,K}(x)) \right]
\\ &\quad\quad\quad \times \Big((\Delta W^{M,K}_m(x))^2-h\sum_{j=1}^{K} \mu_j(\eta_j(x))^2\Big)
\end{split}
\end{equation}
for $m = 0,1,...M-1, \, x \in (0,1)$, the scheme \eqref{eq:concrete_scheme} can be rewritten as
\begin{equation}\label{eq:ImplementationY}
Y_{m+1}^{N,M,K} = P_N \Big(S(h) \, \zeta_m \Big) = \sum_{j=1}^{N} \left\langle e^{A h}\, \zeta_m, e_{j}  \right\rangle_H e_{j} = \sum_{j=1}^{N} e^{-\lambda_{j} \,h} \left\langle \zeta_m, e_{j}  \right\rangle_H e_{j}
\end{equation}
for $m=0,1,...,M-1$ and
$$
Y_0^{N,M,K} = P_N(\xi) = \sum_{j=1}^{N} \left\langle \xi, e_{j}  \right\rangle_H e_{j}.
$$
Here $H= L^2((0,1); \mathbb{R})$, $A$ is the Laplacian with Dirichlet boundary condition times a constant $k>0$ and thus its eigenpairs are given by
$$
e_{j}(x) = \sqrt{2} \sin (j\pi x) \quad \mbox{and} \quad \lambda_{j} = k \pi^2 j^2 \quad \mbox{for}\quad x \in (0,1),\: j\in \mathbb{N}.
$$
For each Fourier mode, we obtain from \eqref{eq:ImplementationY} that
\begin{equation}\label{eq:Implem2}
\left\langle Y_{m+1}^{N,M,K},  e_{j} \right\rangle_H = e^{-\lambda_{j} \,h} \left\langle \zeta_m, e_{j}  \right\rangle_H = \sqrt{2}\, e^{-\lambda_{j} \,h} \int_0^1 \zeta_m(x) \sin (j\pi x) dx, \quad j = 1,2,...,N
\end{equation}
for $m=0,1,...,M-1$ and
\begin{equation}
\left\langle Y_{0}^{N,M,K},  e_{j} \right\rangle_H = \left\langle \xi,  e_{j} \right\rangle_H = \sqrt{2}\, \int_0^1 \xi(x) \sin (j\pi x) dx,  \quad j = 1,2,...,N.
\end{equation}
Therefore, the implementation procedure goes as follows. Given $Y_{m}^{N,M,K}$, one can obtain $\zeta_m$ by \eqref{eq:zeta}. Then we use some numerical integration method (here we choose composite trapezoidal formula) to approximate  $\left\langle \zeta_m, e_{j}  \right\rangle_H$ in \eqref{eq:Implem2} for $j=1,2,...,N$. With $\left\langle \zeta_m, e_{j}  \right\rangle_H, j=1,2,...,N$ at hand, we can get $Y_{m+1}^{N,M,K}$ by \eqref{eq:ImplementationY}. Since the eigenfunction $e_{j}(x) = \sqrt{2} \sin (j\pi x)$ are sine functions, we can invoke built-in functions "idst" and "dst" in matlab to perform efficient computations. Recall that "dst" is a discrete sine transform, which transforms $N$ real numbers $z(k), k=1,2,...,N$ to $N$ real numbers $y(j), j=1,2,...,N$ according to the following formula
\begin{equation}\label{eq:dst}
y(j)= \sum_{k=1}^{N} z(k) \sin(j\pi \frac{k}{N+1}), \quad j=1,...,N,
\end{equation}
and that the "idst" function is an inverse discrete sine transform, which transforms $N$ real numbers $z(k), k=1,2,...,N$ to $N$ real numbers $y(j), j=1,2,...,N$ according to the following formula
\begin{equation}\label{eq:idst}
y(j)= \frac{2}{N+1}\sum_{k=1}^{N} z(k) \sin(j\pi \frac{k}{N+1}), \quad j=1,...,N.
\end{equation}
Setting $z(k) = \zeta_m(\frac{k}{N+1})$ for $k = 1,2,...,N$ in \eqref{eq:idst}, $y(j)$ in \eqref{eq:idst} is in fact a composite trapezoidal formula to numerically approximate $\sqrt{2}\, \left\langle \zeta_m, e_{j}  \right\rangle$. Hence
\begin{equation} \label{eq:appro.integral}
\left\langle Y_{m+1}^{N,M,K},  e_{j} \right\rangle_H \approx e^{-\lambda_{j} \,h} \times y(j)/ \sqrt{2}, \quad j=1,2,...,N.
\end{equation}
To be precise, only $N$ function values of $Y_{m}^{N,M,K}$ at grid points $\frac{k}{N+1}, k = 1,2,...,N$ are used to calculate $N$ function values of $\zeta_m$ at grid points $\frac{k}{N+1}, k = 1,2,...,N$. And then "idst" is used to approximate $\langle Y_{m+1}^{N,M,K},  e_{j}\rangle_H$ by \eqref{eq:appro.integral} for $j=1,2,...,N$.  After that, we use a discrete sine transform "dst" to calculate $N$ function values of $Y_{m+1}^{N,M,K}$ at grid points, which are used to get $N$ function values of $\zeta_m$ at grid points before carrying out numerical integration at next step. Repeating this procedure we can finally obtain $Y_{M}^{N,M,K}$. In Figure \ref{code}, we give the detailed implementation code of our new scheme \eqref{eq:concrete_scheme} for the first test example.

Before closing this subsection, we would like to give some remarks. Here and below, the aliasing errors caused by using composite trapezoidal formula are neglected and are not analyzed mathematically. As illustrated in the following numerical simulations, such errors do not effect the order of convergence.

\begin{figure}[htp]
\centering
\begin{verbatim}
N = 128; M = N^2; A = -pi^2*(1:N).^2/200; mu=(1:N).^-2;
f = @(x) 1-2*x; g = @(x) (x+sin(x).^3)./((1+x.^2).^2);
Y = zeros(1,N);
eta = zeros(1,N); SqrM=sqrt(M);
for n=1:N
    eta = eta + 2*sin(n*(1:N)/(N+1)*pi).^2*mu(n);
end
for m = 1:M
    y = dst(Y)*sqrt(2);
    dW = dst(randn(1,N).*sqrt(mu*2/M));
    g_eva = g(y);
    y = y + f(y)/M+g_eva.*dW + 0.5*SqrM*(g(y+g_eva/SqrM)-g_eva).*(dW.^2-eta/M);
    Y = exp(A/M).*idst(y)/sqrt(2);
end
plot((0:N+1)/(N+1),[0,dst(Y)*sqrt(2),0]);
\end{verbatim}
\caption{Matlab code to simulate one path by the Runge-Kutta type scheme \eqref{eq:concrete_scheme} applied to the SPDE \eqref{eq:spde} with parameters as \eqref{parameter1}.}
\label{code}
\end{figure}

\subsection{Numerical tests}

As the first numerical experiment, we consider an example SPDE \eqref{eq:spde} in the introduction part with initial data $\xi(x) \equiv 0$ for $x \in (0,1)$, $T=1, k=\frac{1}{200}$, and
\begin{equation}\label{parameter1}
\: f(x,y)= 1-2y, \: g(x,y)=\frac{y+\sin^3(y)}{(1+y^2)^2}, \: \mu_j = \frac{1}{j^2}, \: \eta_j(x) = e_j(x) = \sqrt{2}\sin(j \pi x)
\end{equation}
for all $x \in (0,1)$, $y \in \mathbb{R}$, $j \in \mathbb{N}$. Similarly to \cite[Section 4]{JR10b}, one can show that in this case the conditions \eqref{condition_AA} and \eqref{condition_AAA} are fulfilled for $\delta \in (0,\frac{1}{4}), \alpha \in (0, \frac{3}{4}), \gamma \in (\frac{1}{2}, \frac{3}{4})$. As a result, Assumption \ref{semigroup}-\ref{BB_assum} are all satisfied for $\beta = \frac{1}{5}, \delta \in (0,\frac{1}{4}), \alpha \in (0, \frac{3}{4}), \gamma \in (\frac{1}{2}, \frac{3}{4})$.

According to the computational analysis in \cite{JR10b}, we know that the linear implicit Euler method \eqref{eq:Escheme} with $M=N^3, K=N$ promises the existence of some real constants $C_r >0$ and arbitrarily small $r \in (0, \frac{3}{2})$ such that
\begin{equation}\label{eq:Euler_error}
\left(\mathbb{E} \left\|X_T-\bar{Y}_{N^3}^{N,N^3,N}\right\|_H^2 \right)^{\frac{1}{2}} \leq C_r \, N^{r-\frac{3}{2}}.
\end{equation}

For the infinite version of Milstein type method \eqref{eq:CMscheme}, it is shown in \cite{JR10b} that, $N^2$ time steps, in contrast to $N^3$ time steps for the linear implicit Euler scheme \eqref{eq:Escheme}, are required to achieve \eqref{eq:Euler_error}, that is, \eqref{eq:CMscheme} with $N^2$ time steps guarantees that for some real
constants $C_r >0$ and arbitrarily small $r \in (0, \frac{3}{2})$
\begin{equation}
\left(\mathbb{E} \Big\|X_T-\tilde{Y}_{N^2}^{N,N^2,N}\Big\|_H^2 \right)^{\frac{1}{2}} \leq C_r \, N^{r-\frac{3}{2}}.
\end{equation}

For the new scheme \eqref{eq:concrete_scheme} applied to \eqref{eq:spde}, the main result (Theorem \ref{thm:mainresult}) in this paper shows that
for some positive constants $C_r$
\begin{equation}\label{eq:convergence}
    \left(
      \mathbb{E}
      \left\|
        X_{T}
        - Y^{N,M,K}_M
      \right\|_H^2
    \right)^\frac{1}{2}
    \leq
    C_r
  \left(
      N^{r-\frac{3}{2}}
    + K^{r-\frac{3}{2}}
    +
      M^{r-\frac{3}{4}}
  \right)
\end{equation}
holds for all arbitrarily small $r \in (0, \frac{3}{4})$.
Choosing $M=N^2, K=N$ in the preceding result gives
\begin{equation}\label{eq:RK_error}
\left(\mathbb{E} \Big\|X_T-Y_{N^2}^{N,N^2,N}\Big\|_H^2\right)^{\frac{1}{2}} \leq C_r \, N^{r-\frac{3}{2}}
\end{equation}
for arbitrarily small $0< r < \frac{3}{2}$.

In Figure \ref{code}, we present detailed implementation code of the scheme \eqref{eq:concrete_scheme}. The term $h\sum_{j=1}^{K} \mu_j(\eta_j)^2$ in \eqref{eq:concrete_scheme} is computed once in advance for which $O(N^2)$ computational operations are needed. After that,
$O(N\log(N))$ further computational operations and independent standard normal random variables are needed to compute $Y_{m+1}^{N,N^2,N}$ from $Y_{m}^{N,N^2,N}$ by using the fast Fourier transform. Since $N^2$ time steps are used, $O(N^3\log(N))$ computational operations and random variables are required to obtain $Y_{N^2}^{N,N^2,N}$. Taking the convergence order $\frac{3}{2}-$ in \eqref{eq:RK_error} into account shows that the scheme \eqref{eq:concrete_scheme}
promises the overall convergence order $\frac{1}{2}-$, which is the same as that of Milstein type scheme \eqref{eq:CMscheme}. But for the linear implicit Euler scheme \eqref{eq:Escheme}, $N^3$ time steps are used and one can just get an overall convergence order of $\frac{3}{8}-$.

\begin{figure}[htp]
         \centering
         \includegraphics[width=4in,height=3in]{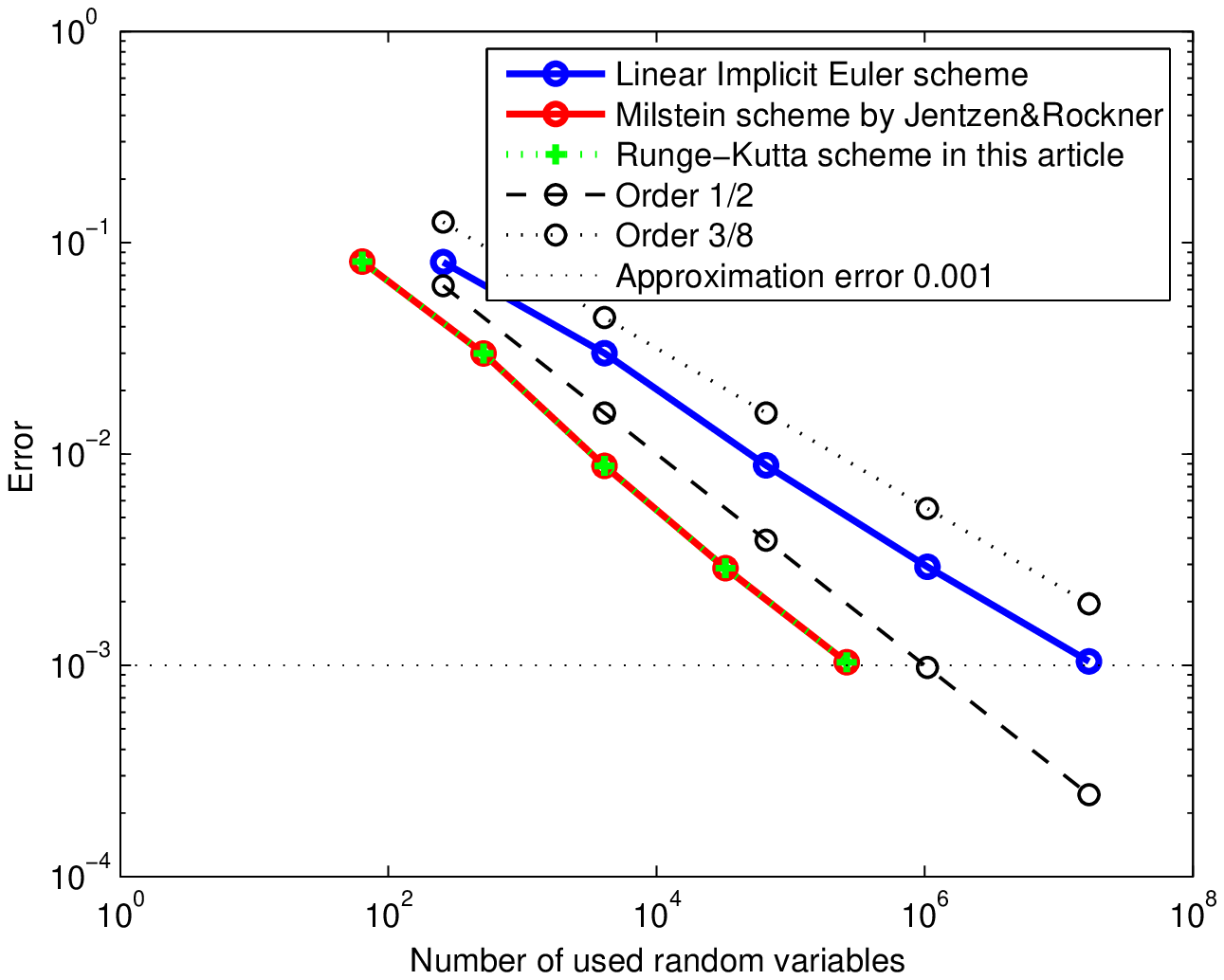}
         \caption{Numerical results for SPDE \eqref{eq:spde} with parameters as \eqref{parameter1}.}
         \label{fig1}
\end{figure}
Figure \ref{fig1} depicts approximation errors \eqref{eq:strong_error} of the various approximations $\bar{Y}_{N^3}^{N,N^3,N}$, $\tilde{Y}_{N^2}^{N,N^2,N}, Y_{N^2}^{N,N^2,N}$ with $N =4,8,16,32,64$ against the number of used normal random variables on a log-log scale. As a measure for the computational effort, here we take the number of realizations of independent random variables needed for the calculation of the approximation. One can detect that the numerical results are consistent with our assertions on the convergence order. Besides, the Runge-Kutta method \eqref{eq:concrete_scheme} and the Milstein type scheme \eqref{eq:CMscheme} produce nearly the same approximation error. Numerical results also show that both the Runge-Kutta method \eqref{eq:concrete_scheme} and the Milstein type scheme \eqref{eq:CMscheme} are much more computationally effective than the linear implicit Euler scheme \eqref{eq:Escheme}. For instance,  $\bar{Y}_{64^3}^{64,64^3,64}$ in the case of the linear implicit Euler scheme and  $\tilde{Y}_{64^2}^{64,64^2,64}, Y_{64^2}^{64,64^2,64}$ in the case of the Runge-Kutta method and the Milstein type scheme achieve the precision $\varepsilon = 0.001$ in \eqref{eq:strong_error}. For one path, it needs to generate $64^4 = 16777216$ independent normal random variables and costs 105.312000 seconds to simulate $\bar{Y}_{64^3}^{64,64^3,64}$. But the number of random variables needed to generate decreases to $64^3=262144$ for one path simulation of $\tilde{Y}_{64^2}^{64,64^2,64}, Y_{64^2}^{64,64^2,64}$. Accordingly, the runtime for one path simulation of $\tilde{Y}_{64^2}^{64,64^2,64}$ and $Y_{64^2}^{64,64^2,64}$ are, respectively, reduced to 2.328000  seconds and 1.984000 seconds. In Table \ref{table1}, we list runtime (seconds) for one path simulation using the three methods with various $N$ ($N = 32,64,128,256$). Note that the "exact" mild solution is identified with the numerical solution using very small stepsize and that the matlab codes to simulate $\bar{Y}_{N^3}^{N,N^3,N}, \tilde{Y}_{N^2}^{N,N^2,N}$ and $Y_{N^2}^{N,N^2,N}$ are presented in Figure 2, Figure 3 from \cite{JR10b} and Figure \ref{code} in this article, respectively. It turns out that the schemes \eqref{eq:concrete_scheme} and \eqref{eq:CMscheme} are progressively faster than the linear implicit Euler method \eqref{eq:Escheme} as $N$ increases. Also, we observe from Table \ref{table1} that the Runge-Kutta type scheme \eqref{eq:concrete_scheme} is faster than the Milstein type scheme \eqref{eq:CMscheme}. This is due to the fact that evaluation of the partial derivative $(\frac{\partial }{\partial y}g)(x,y)=\frac{3(1+y^2)\sin^2(y)\cos(y)-4y\sin^3(y)-3y^2+1}{(1+y^2)^3}$ costs more time than evaluation of the function $g(x,y)=\frac{y+\sin^3(y)}{(1+y^2)^2}$.

\begin{table}[htp]
\begin{center} \footnotesize
\caption{Runtime (seconds) for one path simulation using the three schemes $\bar{Y}_{m}^{N,N^3,N}, \tilde{Y}_{m}^{N,N^2,N}, Y_{m}^{N,N^2,N}$ with $N=32,64,128,256$ }
\label{table1}
\begin{tabular*}{16cm}{@{\extracolsep{\fill}}cccc}
\hline\\
 & Linear Implicit Euler scheme  & Milstein scheme  &  Runge-Kutta scheme  \\
\hline\\
$N=32$ & 10.172000 &  0.547000 &  0.469000  \\
$N=64$ & 105.312000 & 2.328000 & 1.984000  \\
$N=128$ & 1068.969000  &12.547000 & 10.344000  \\
$N=256$ & 12604.844000 &77.016000 &60.312000  \\
\hline
\end{tabular*}
\end{center}
\end{table}
As the second numerical experiment, we consider the case,
where the two operators $A$ and $Q$ do not share the same eigenfunctions.
More accurately, we choose the initial data $\xi(x) \equiv 0$ for $x \in (0,1)$, $T=1, k=\frac{1}{50}$, and
the other parameters are set as
\begin{equation}\label{parameter2}
f(x,y)= 1-y, \: g(x,y)=\frac{y}{1+y^2}, \: \mu_0=0, \: \mu_j = \frac{1}{j^3}, \: \eta_0(x)=1,\: \eta_j(x) = \sqrt{2}\cos(j \pi x)
\end{equation}
for all $x \in (0,1)$, $y \in \mathbb{R}$, $j \in \mathbb{N}$.

\begin{figure}[htp]
         \centering
         \includegraphics[width=4in,height=3in]{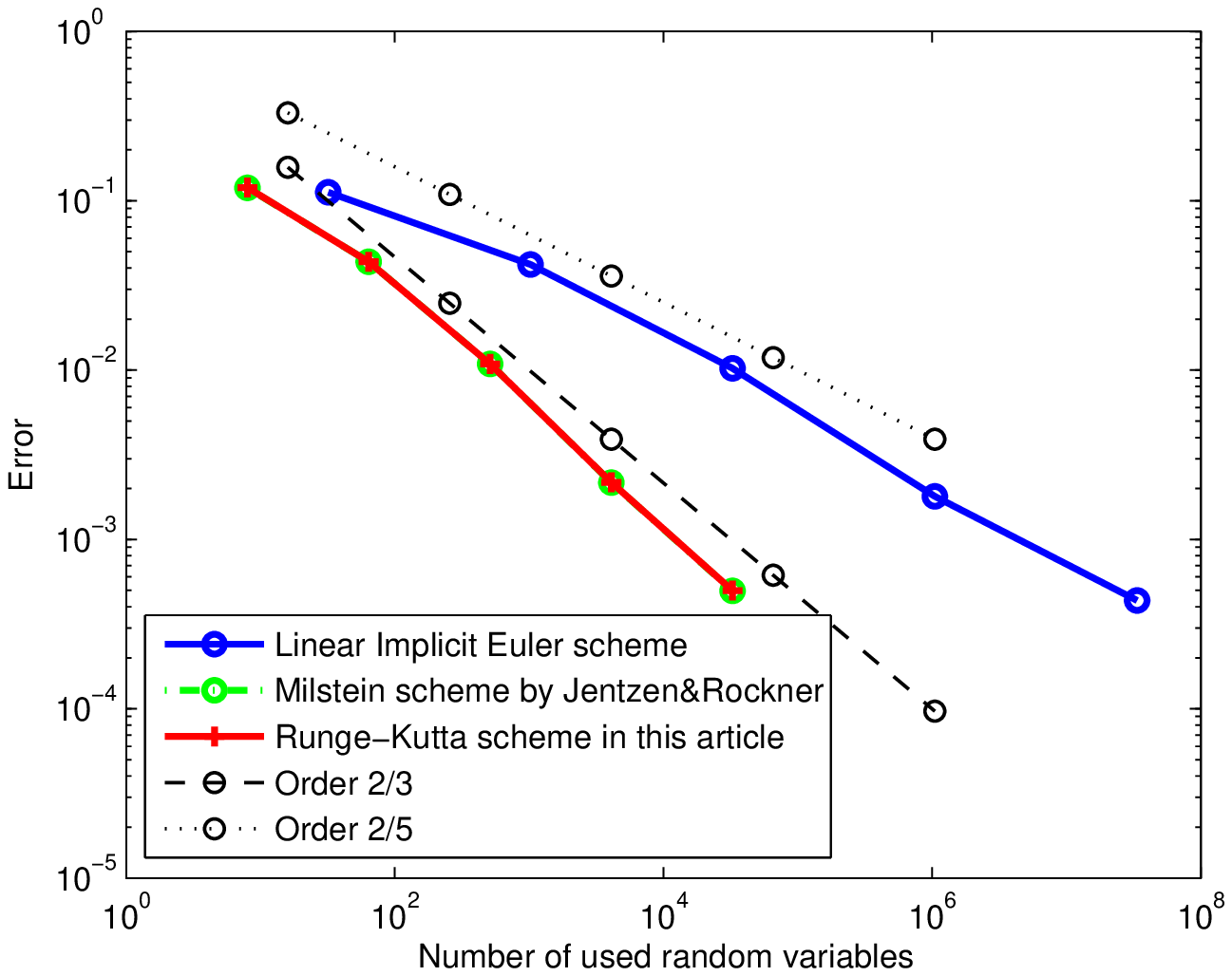}
         \caption{Numerical results for SPDE \eqref{eq:spde} with parameters as \eqref{parameter2}.}
         \label{fig3}
\end{figure}

For this example, it is shown in \cite{JR10b} that Assumption \ref{semigroup}-\ref{BB_assum} are all satisfied with $\beta = \frac{1}{5}, \delta \in (0,\frac{1}{2}), \alpha \in (0, \frac{2}{3}), \gamma \in (\frac{1}{2}, 1)$. Consequently, in this case Theorem \ref{thm:mainresult} shows that
\begin{equation}\label{eq:convergence2}
    \left(
      \mathbb{E} \Big\|
        X_{T}
        - Y^{N,N^2,N}_{N^2}
      \Big\|_H^2
    \right)^\frac{1}{2}
    \leq
    C_r \,
      N^{r-2}
\end{equation}
holds for some positive constants $C_r$, all arbitrarily small $r \in (0, 2)$
and $N \in \mathbb{N}$. Hence its overall convergence order is $\frac{2}{3}-$. For the linear implicit Euler scheme \eqref{eq:Escheme},
\begin{equation}\label{eq:Euler_convergence2}
    \left(
      \mathbb{E} \Big\|
        X_{T}(x)
        - \bar{Y}^{N,N^4,N}_{N^4}(x)
      \Big\|_H^2
    \right)^\frac{1}{2}
    \leq
    C_r \,
      N^{r-2}
\end{equation}
holds for some positive constant $C_r$, all arbitrarily small $r \in (0, 2)$
and $N \in \mathbb{N}$. \eqref{eq:Euler_convergence2} implies that the
linear implicit Euler scheme has the overall convergence order of
$\frac{2}{5}-$.
These asymptotic results can be observed clearly in Figure \ref{fig3}, where
approximation errors  of the three approximations
$\bar{Y}_{N^4}^{N,N^4,N}, \tilde{Y}_{N^2}^{N,N^2,N}, Y_{N^2}^{N,N^2,N}$ with
$N =2,4,8,16,32$ against the number of used normal random variables are plotted.

%

%
%
%
%
\bibliographystyle{acm}

\end{document}